\documentclass[twocolumn,letter]{autart}    

\usepackage{amsmath,amssymb,color}
\usepackage{graphicx}
\usepackage{wrapfig} 

\usepackage{latexsym}
\usepackage{verbatim}
\usepackage{cite}
\usepackage[latin1]{inputenc}
\usepackage{enumerate}
\usepackage{enumitem}
\renewcommand{\theenumi}{(\roman{enumi})}

\usepackage[font=footnotesize]{caption}
\usepackage{arydshln}

\usepackage{amsfonts}
\usepackage{amsbsy}
\usepackage{bm}
\usepackage{bbm}
\usepackage{pmat}
\usepackage[numbers]{natbib}

\newcommand{\longthmtitle}[1]{\mbox{}{\rm \textbf{(#1)}}}
\newcommand\oprocendsymbol{\hbox{$\bullet$}}
\newcommand\oprocend{\relax\ifmmode\else\unskip\hfill\fi\oprocendsymbol}

\newcommand{\QED}{\hspace*{\fill} $\blacksquare$}

\DeclareMathOperator{\im}{im}

\DeclareMathOperator{\conv}{conv}

\DeclareMathOperator{\osL}{osLip}

\DeclareMathOperator*{\argmin}{arg\,min}

\newcommand{\norm}[1]{\ensuremath{\| #1 \|}}

\newcommand{\until}[1]{\{1,\dots,#1\}}

\DeclareMathOperator{\interior}{int}
\DeclareMathOperator{\uniform}{Uni}

\let\leq\leqslant
\let\geq\geqslant

\newcommand{\R}{\mathbb R}
\newcommand{\N}{\mathbb N}

\newcommand{\lse}{\operatorname{lse}}

\newcommand{\calF}{\ensuremath{\mathcal{F}}}
\newcommand{\calG}{\ensuremath{\mathcal{G}}}

\newcommand{\calS}{\ensuremath{\mathcal{S}}}

\newcommand{\calZ}{\ensuremath{\mathcal{Z}}}

\newcounter{todocounter}
\usepackage[colorinlistoftodos]{todonotes}

\newtheorem{theorem}{Theorem}[section]
\newtheorem{proposition}[theorem]{Proposition}
\newtheorem{lemma}[theorem]{Lemma}
\newtheorem{corollary}[theorem]{Corollary}

\newtheorem{remark}[theorem]{Remark}

\newtheorem{assumption}{Assumption}

\newtheorem{example}{Example}
\newtheorem{problem}{Problem}
\newtheorem{algo}{Algorithm}

\begin{document}




\begin{frontmatter}
  \title{Cautious optimization via data informativity}

  \thanks[footnoteinfo]{A preliminary version of this work was
    submitted as~\cite{JE-JC:23-cdc} to the IEEE Conference on
    Decision and Control. Corresponding author J.~Eising.}
  
  \author[First]{Jaap Eising}
  \quad 
  \author[Second]{Jorge Cort\'{e}s}            
  
  \address[First]{Automatic Control Laboratory, ETH Z\"{u}rich,
    jeising@ethz.ch}
  \address[Second]{Department of Mechanical and Aerospace
    	Engineering, University of California, San Diego,
    	cortes@ucsd.edu \vspace{-2em}}
  \begin{abstract}
    This paper deals with the problem of accurately determining
    guaranteed suboptimal values of an unknown cost function on the
    basis of noisy measurements.  We consider a set-valued variant to
    regression where, instead of finding a best estimate of the cost
    function, we reason over all functions compatible with the
    measurements and apply robust methods explicitly in terms of the
    data.  Our treatment provides data-based conditions under which
    closed-forms expressions of upper bounds of the unknown function
    can be obtained, and regularity properties like convexity and
    Lipschitzness can be established.  These results allow us to
    provide tests for point- and set-wise verification of
    suboptimality, and tackle the cautious optimization of the unknown
    function in both one-shot and online scenarios.  We showcase the
    versatility of the proposed methods in two control-relevant
    problems: data-driven contraction analysis of unknown nonlinear
    systems and suboptimal regulation with unknown dynamics and cost.
    Simulations illustrate our results.
  \end{abstract}
\end{frontmatter}

\section{Introduction}

In many real-world applications involving optimization, the cost or
reward function is not fully known and has to be ascertained from
measurements. Such situations can arise due to a variety of factors,
including system complexity, large-scale structure, or lack of access
to the relevant information in unknown or adversarial scenarios.
These considerations motivate the need for quantifiable performance
guarantees that can be used for control and optimization in
safety-critical applications. In this paper, we are therefore
interested in cautious optimization of unknown functions, i.e.,
accurately determining guaranteed suboptimal values on the basis of
noisy measurements. We consider both one-shot scenarios, where we are
required to make decisions on the basis of a given set of data, and
online scenarios, where the suboptimization of the unknown function
can be repeatedly refined based on the collection of additional
measurements. Our approach, based on data informativity, allows us to
derive results which are robust against worst-case situations.

\emph{Literature review:} Data-based optimization is a widely
applicable and thoroughly investigated area. Giving a full overview of
all methods and results is not possible, and thus we focus here on
those most relevant to our
treatment. 
The problems considered here are robust optimization programs, see
e.g.~\cite{AB-LEG-AN:09,DB-DBB-CC:11} and references therein.
Solution approaches to data-based optimization problems usually start
by using the data to obtain an estimate of the unknown function which
in some sense `best' explains the measurements. This requires the
definition of the class of admissible functions. Popular techniques
are based on Gaussian processes, as in~\cite{CER-CKIW:05} or the
approximation properties of neural networks~\cite{GC:89,PT-BG:23}.
{\color{black} Of particular relevance to our analysis here is the work
  on set-membership estimation
  \cite{MM-AV:91,MM-CN:02,JC-SJR-CER-JM:20}, which, instead of
  characterizing the best function, considers model-free set-valued
  uncertainties using Lipschitz interpolation.

  In contrast to these approaches, we consider the situation where the
  unknown function can be written as a linear combination of a set of
  basis functions. In practice, this is not always possible, and thus
  the choice of basis induces a model mismatch error, which needs to
  be taken into account. Common choices that balance expressiveness
  and the number of parameters are linear, polynomials, Gaussian, and
  sigmoidal functions.  Depending on the problem, one may also choose
  a basis by taking into account physical considerations or employing
  spectral methods~\cite{MB-RM-IM:12,MK-IM:20,MH-JC:23-auto}. In the
  specific case of nonlinear system identification, an overview of the
  problem of selecting a suitable parametrization can be found
  in~\cite{ON:20}.  Once a basis is chosen, the data is used to
  perform regression~\cite{TH-RT-JF:13,RT:96} on the parameters. What
  remains to be defined then is the metric chosen to decide the best
  estimate. Common methods are least squares (minimal Frobenius norm),
  ridge regression (minimal $L_2$ norm), and sparse or Lasso
  regression (minimal $L_1$ norm). The nonlinear system-identification
  literature has also brought forward methods which determine models
  that are sparse~\cite{SLB-JLP-JNK:16}, of low rank (via dynamic mode
  decomposition)~\cite{PJS:10,PJS:22,JNK-SLB-BWB-JLP:16}, or
  both~\cite{MRJ-PJS-JWN:14}. Regardless of how one obtains an
  estimate, the second step in data-based optimization is to employ
  either a certainty-equivalent or robust method to find optimal
  values.

  The alternative approach we use here is in line with the literature
  of data-driven control methods on the basis of Willems' fundamental
  lemma~\cite{JCW-PR-IM-BLMDM:05}. More specifically, we consider the
  concept of data
  informativity~\cite{HJVW-JE-HLT-MKC:20,HJVW-JE-HLT-MKC:23}. Simply
  put, these methods take the viewpoint that guaranteed conclusions
  can be drawn from data only if \emph{all} systems compatible with
  the measurements have this property. Of specific interest to this
  paper are works that deal with nonlinear systems and a number of
  recent works have worked with systems given in terms of
  bilinear~\cite{IM:22} or polynomial basis
  functions~\cite{MG-CDP-PT:21,CDP-MR-PT:23}, albeit in discrete
  time. Simultaneously, there have been efforts in generalizing
  Willems' lemma to continuous time~\cite{VGL-MAM:22,PR-MKC-HJVW:22}
  and linking the results in discrete and continuous time by
  investigating sampling~\cite{JE-JC:23-tac}. A related approach is to
  perform system identification while keeping in mind the uncertainty
  bounds used in (robust) control, which is known as identification
  for control~\cite{MG:05}. The convergence over time of similar
  set-based identification schemes is investigated
  in~\cite{YL-JY-LC-TK-AW:23}.
	
Apart from considering optimization, we apply our results to
data-based contraction analysis of nonlinear
systems~\cite{WL-JJES:98,FB:23} and the regulation of an unknown
system to a suboptimal point of an unknown cost function.}
Data-based regulation with known cost functions has been investigated
in a number of forms, including linear quadratic regulation,
see~\cite{FD-PT-CDP:23} and references therein, and more general cost
functions via optimization-based
controllers~\cite{GB-MV-JC-EDA:24-tac,AH-ZH-SB-GH-FD:24}. In~\cite{LC-GB-ED:23},
such controllers have been employed to optimize unknown cost
functions.  Finally, the online optimization scenario considered here
is inspired by methods such as extremum-seeking
control~\cite{MK-HHW:00,KBA-MK:03,ART-DP:01}. However, our
implementation and performance guarantees are markedly different.

\emph{Statement of contributions:} We tackle the data-based
optimization of an unknown vector-valued function with the goal of
obtaining quantifiable performance guarantees.  Starting from the
assumption that the unknown function is parameterized in terms of a
finite number of basis functions, we consider \emph{set-valued
  regression}, i.e., we reason with the set of \textit{all} parameters
compatible with the measurements for which the noise satisfies a bound
given in terms of a quadratic matrix inequality.  We take a worst-case
approach regarding noise realizations, leading us to provide solutions
to the following \textit{cautious optimization} problems:
\begin{enumerate}
\item We provide conditions in terms of the measurements, basis
  functions, and noise model which guarantee that either the $2-$norm
  or a linear functional of the unknown function is upper bounded by a
  given $\delta\in\mathbb{R}$ for a given
  $z\in\mathbb{R}^n$. Moreover, these results make it possible for us
  to find the smallest worst-case bounds;
\item We also provide data-based conditions for convexity of either
  the true function or its worst-case norm bounds. These allow us to
  introduce a gradient-based method for finding the minimizer of the
  upper bound;
\item We also identify data-based conditions for Lipschitzness which,
  using interpolation, allows us to derive upper bounds guaranteed
  to hold over compact sets.
\end{enumerate}
We illustrate the versatility of the results in two application
scenarios: contraction analysis of nonlinear systems (for which we
provide data-based tests to determine one-sided Lipschitz constants
and use them to guarantee contractivity) and suboptimal regulation of
unknown plants with unknown cost functions (for which we employ
convexity and Lipschitzness to find a fixed input such that the value
of the cost function at the corresponding fixed point can be
explicitly bounded above).  Our final contribution is to online
data-based optimization.  We describe a framework for incorporating
new measurements of the unknown function which includes online local
descent, a method that iteratively collects measurements locally and
then minimizes an upper bound of the form described above.  We prove
that under mild assumptions on the collection of data, the set of
parameters consistent with all collected measurements shrinks. This
allows us to conclude that the upper bounds found with online local
descent converge to the true optimizer of the unknown function.

Preliminary results of this paper were presented in the conference
article~\cite{JE-JC:23-cdc}, whose focus was restricted to scalar
functions and provided upper bounds and analyzed convexity on the
basis of data. The conference paper also presented a simplified
version of online local descent for the scalar case.  All of these are
special cases of the present work. The generalization here from scalar
to vector-valued functions is instrumental in expanding the range of
possibilities for application in analysis and control of the proposed
data-based optimization framework.

\emph{Notation:} Throughout the paper, we use the following
notation. We denote by $\N$ and $\R$ the sets of nonnegative integers
and real numbers, respectively. We let $\R^{n\times m}$ denote the
space of $n\times m$ real matrices and $\mathbb{S}^n$ the space of
symmetric $n\times n$ matrices. For a vector $x\in\mathbb{R}^n$ we
denote the standard $p-$norm by $\norm{x}_p$. Given {\color{black} an
  invertible matrix} $R\in\mathbb{R}^{n\times n}$, we define the
weighted $p-$norm by $\norm{x}_{p,R}:=\norm{Rx}_p$. The
$i^{\textrm{th}}$ standard basis vector {\color{black} of
  $\mathbb{R}^n$} is denoted $e_i$. The induced matrix norm w.r.t. to
the $p-$norm is denoted $\norm{A}_p$ for $A\in\mathbb{R}^{n\times
  m}$. For vectors $v\in\R^n$, we write $v\geq0$ (resp. $v>0$) for
elementwise nonnegativity (resp. positivity). The sets of such vectors
are denoted $\mathbb{R}^n_{\geq0}:=\{v\in\R^n | v\geq0\}$ and
$\mathbb{R}^n_{>0}:=\{v\in\R^n | v>0 \}$. For $P\in\mathbb{S}^{n}$,
$P\geq 0$ (resp. $P> 0$) denotes that $P$ is symmetric positive
semi-definite (resp. definite). For a {\color{black} symmetric} matrix
$M$, we denote $\lambda_{\textup{min}}(M)$ and
$\lambda_{\textup{max}}(M)$ for the smallest and largest eigenvalue of
$M$. We denote the smallest singular value of
$M\in\mathbb{R}^{n\times m}$ by $\sigma_-(M)$ and its Moore-Penrose
pseudoinverse by $M^\dagger$.  The convex hull and interior of
$\calS\subseteq\mathbb{R}^n$ are denoted by $\conv(\calS)$ and
$\interior(\calS)$, resp. Given a symmetric matrix
$M\in\mathbb{S}^{u+v}$, when $u$ and $v$ are clear from context, we
partition it as:
\vspace{-0.5em}\[
  M = {\scriptsize	
  \begin{bmatrix}
    M_{11} & M_{12}
    \\
    M_{21} & M_{22}
  \end{bmatrix}},\vspace{-0.5em}
\] 
where $M_{11}\in\mathbb{S}^{u}$ and $M_{22}\in
\mathbb{S}^{v}$. {\color{black} If $M_{22}\leq 0$ and
  $\ker M_{22} \subseteq \ker M_{12}$, we denote the (generalized)
  Schur complement of $M$ with respect to $M_{22}$ by
  $M|M_{22}:= M_{11}- M_{12} M_{22}^\dagger M_{21}$.}  Lastly,
{\color{black} such} a partitioned matrix gives rise to a quadratic
matrix inequality (QMI), whose set of solutions is denoted
\vspace{-0.6em}\[\vspace{-1.2em}
  \calZ(M) := \left\{ Z\in \mathbb{R}^{v\times u} \mid
  {\scriptsize	\begin{bmatrix}
    I_u
    \\
    Z
  \end{bmatrix}^\top}
\!\!
M
{\scriptsize\begin{bmatrix}
  I_u
  \\
  Z
\end{bmatrix}}
\geq 0 \right\} .
\]
\vspace{-0.8em}
\section{Problem formulation}\label{sec:problem}
Consider a collection of known \emph{basis functions} (or
\emph{features}), denoted $\phi_i:\mathbb{R}^n\rightarrow\mathbb{R}$
for $i=1,\ldots, k$. We collect the basis
functions into a vector-function as
\vspace{-0.3em}
\begin{align}
  \label{eq:b}
  \vspace{-0.3em}  b:\mathbb{R}^n\rightarrow\mathbb{R}^k, \quad b(z) :=
  {\scriptsize	\begin{bmatrix}
    \phi_1(z)^\top & \cdots & \phi_k(z)^\top 
  \end{bmatrix}^\top}.  
\end{align}
{\color{black} Using these basis functions}, we define linear combinations
as parameterized by $\theta\in \mathbb{R}^{k\times m}$:
\vspace{-0.3em}
\[
  \vspace{-0.3em}
  \phi^\theta: \mathbb{R}^n\rightarrow\mathbb{R}^m, \quad
  \phi^\theta(z)= \theta^\top b(z).
\]
Our starting point is an unknown function
$\hat{\phi}:\mathbb{R}^n\rightarrow\mathbb{R}^m$ that can be expressed
in this form: i.e., there exists
$\hat{\theta}\in\mathbb{R}^{k\times m}$ such that
$\hat{\phi}(z) = \hat{\theta}^\top b(z)$.

Our goal is to investigate properties and {\color{black} find optimal
  values of} the function $\hat{\phi}$ on the basis of
measurements at points $\{z_i\}_{i=1}^T$. We assume the measurements
are \emph{noisy}, that is, we collect $y_i =\hat{\phi}(z_i) + w_i$,
where $w_i$ denotes an unknown disturbance vector for each~$i$. To
express this in compact form, define the matrices
$Y,W\in \mathbb{R}^{m\times T}$ and $\Phi \in \mathbb{R}^{k\times T}$
by
\vspace{-0.5em}\begin{equation}\label{eq:Y-W}
  \begin{split}
  & Y := \lbrack\begin{matrix} y_1 & \cdots &
    y_T \end{matrix}\rbrack, \; \quad  W := \lbrack\begin{matrix} w_1 & \cdots &
    w_T \end{matrix}\rbrack, \\
 &\Phi := \lbrack\begin{matrix} b(z_1) & \cdots & b(z_T)\end{matrix}\rbrack
    = {\scriptsize	 \begin{bmatrix} \phi_1(z_1) & \dots & \phi_1(z_T) \\ \vdots & &
    	\vdots \\ \phi_k(z_1) & \dots & \phi_k(z_T) \end{bmatrix}}.
\end{split}
\end{equation}
\vspace{-1.5em}

Then, for the true value of $\hat{\theta}$, we have 
$Y = \hat{\theta}^\top\Phi +W$. In this equation, the matrices $Y$ and 
$\Phi$ are known, and $\hat{\theta}$ and $W$ are unknown.

A common line of reasoning (see e.g.,~\cite{TH-RT-JF:13}) to
approximate the unknown function $\hat \phi$ is the
following. {\color{black}Assuming that small noise samples are more
  probable than large noise samples, we attempt to find a solution
  $\theta$ to $Y = \theta^\top\Phi +W$ for which the Frobenius norm of
  $W$, denoted $\norm{W}_F$, is as small as possible. This value is
  attained for $\theta= (Y\Phi^\dagger)^\top$.}

Here instead we consider bounded noise and reason with the set of
functions consistent with the measurements. Formally, we assume that
the noise conforms to a model defined by a QMI
given in terms of a partitioned matrix. The following assumption 
describes the noise model.

\begin{assumption}[Noise model]\label{as:noise model}
  The noise samples satisfy $W^\top\in \calZ(\Pi)$,
  where $\Pi\in\mathbb{S}^{m+T}$ is such that 
  $\Pi_{22}<0$ and $\Pi| \Pi_{22}\geq 0$.
\end{assumption}

According to~\cite[Thm. 3.2]{HJVW-MKC-JE-HLT:22}, under
Assumption~\ref{as:noise model}, the set $\calZ(\Pi)$ is nonempty,
convex, and bounded. A common example of such noise models is the case
where $W W^\top = \sum_{i=1}^T w_iw_i^\top \leq Q$ for some $Q\geq 0$.
 {\color{black} This choice is analogous to assuming that $w$ has bounded energy.
 If we, moreover, take $Q= \gamma T I_n$ or $Q = \gamma YY^\top$ we obtain bounds in 
 terms of the horizon, or signal-to-noise bounds. 
 Lastly, as shown 
  in~\cite[Sec. 5.4]{HJVW-MKC-JE-HLT:22}, these noise models can be
  employed as confidence intervals corresponding to a given
  probability in the setting of Gaussian noise}.

Under Assumption~\ref{as:noise model}, one can describe the set of
parameters consistent with the measurements as
\vspace{-0.3em}\begin{equation}\label{eq:def Theta}
\vspace{-0.3em}  \Theta = \{ \theta \in
  \mathbb{R}^{k\times m}
  \mid Y =
  \theta^\top\Phi +W, W^\top\in\calZ(\Pi) \}.
\end{equation}
Note that we can write
\vspace{-0.7em}\[ \vspace{-0.6em} Y=\theta^\top\Phi +W
  \iff {\scriptsize	\begin{bmatrix} I_m & W \end{bmatrix}}= {\scriptsize	\begin{bmatrix} I_m &
    \theta^\top \end{bmatrix} \begin{bmatrix} I_m & Y \\ 0&
    -\Phi \end{bmatrix}}.
\]
Thus, if we define
\vspace{-0.9em}\begin{equation}\label{eq:defN} 
\vspace{-0.3em}  N:=
  {\scriptsize	\begin{bmatrix} I_m & Y \\ 0& -\Phi \end{bmatrix}}
  \Pi {\scriptsize\begin{bmatrix} I_m & Y \\ 0& -\Phi \end{bmatrix}^\top},  
\end{equation} 
it follows immediately that $\Theta = \calZ(N)$. We refer to the
procedure of obtaining the set of parameters $\Theta$ from the 
measurements as \emph{set-valued regression}. 

{\color{black}
  \begin{remark}\longthmtitle{Richness of the
      measurements}\label{rem:exciting}
    {\rm Note that the set $\Theta$ is compact if and only if
      $N_{22}<0$. Since $N_{22} = \Phi\Pi_{22}\Phi^\top$ and
      $\Pi_{22}<0$, this holds if and only if $\Phi$ has full row
      rank. In turn, this requires that the basis functions are not
      identical and that the set of points $\{z_i\}_{i=1}^T$ is `rich'
      enough to distinguish them. If the unknown function represents
      the state equation of a dynamical system, such rank properties
      are often referred to as (persistent) `excitation',
      cf.~\cite{JCW-PR-IM-BLMDM:05}. \oprocend }
\end{remark}}

\begin{remark}\longthmtitle{Least-squares estimate}\label{rem:lse}
  {\rm One can check that, under Assumption~\ref{as:noise model},
    $\theta^{\lse}:=-N_{22}^\dagger N_{21} \in \Theta$. This means
    that the function $ (\theta^{\lse})^\top b(z)$ is consistent with
    the measurements. In fact,
    \vspace{-1em}\[
     \scriptsize N|N_{22}=\begin{bmatrix} I_m \\ -N_{22}^\dagger
        N_{21} \end{bmatrix}^\top
      N \begin{bmatrix} I_m \\
        -N_{22}^\dagger N_{21} \end{bmatrix}\geq \begin{bmatrix} I_m \\
        \theta\end{bmatrix}^\top N \begin{bmatrix} I_m \\
        \theta \end{bmatrix},
    \]
    for any $\theta^\top\in \calZ(N)$. Therefore, $\theta^{\lse}$ is
    the value that maximizes the value of the quadratic
    inequality. This leads us to refer to the function
    $\phi^{\lse}(z;\Theta):=(\theta^{\lse})^\top b(z) =
    -N_{12}N_{22}^\dagger b(z)$ as the \textit{least-squares estimate}
    of $\hat{\phi}(z)$. \oprocend }
\end{remark}

{\color{black} Note that, without further assumptions on the data, we
  cannot distinguish $\hat{\phi}$ from any of the other functions
  $\phi^\theta$ with $\theta \in \Theta$. Hence, we can \textit{only}
  conclude that a given optimality criterion of $\hat{\phi}$ holds if
  it holds for all functions $\phi^\theta$ with $\theta \in
  \Theta$. Hence, a reasonable approach would be to optimize some
  criterion for all $\phi^\theta$ with $\theta \in \Theta$
  \textit{simultaneously}. However, changes in the parameter $\theta$
  might lead to changes in the quantitative behavior and in particular
  the location of optimal values of the corresponding
  function~$\phi^\theta$. This means that a small error in estimating
  the true parameter corresponding to the unknown
  function~$\hat{\phi}$ might lead to significant error. This
  motivates us to focus on suboptimization problems
  instead.} To be precise, we investigate whether, for a given
$\delta \in \mathbb{R}$ and $c\in\mathbb{R}^m$, we can ensure that
\vspace{-0.3em}\begin{equation}\label{eq:problems}
  \vspace{-0.3em}\norm{\hat{\phi}(z)}_2 \leq \delta
  \quad \textrm{ or }\quad c^\top \hat{\phi}(z) \leq
  \delta.
\end{equation}
We formalize this next.

{\color{black}\begin{problem}[Cautious
    optimization]\label{prob:cautious-opt} {\rm Consider an unknown
      function $\hat{\phi}:\mathbb{R}^n\rightarrow\mathbb{R}^m$, a
      noise model $\Pi$ such that Assumption~\ref{as:noise model}
      holds, and data $(Y,\Phi)$ of the true function:
    \begin{enumerate}
    \item\label{prob:verif} (Verification of suboptimality) given
      $z\in\mathbb{R}^n$, determine the minimal value of
      $\delta\in\mathbb{R}$ as a function of the data $(Y,\Phi)$, for
      which we can guarantees that either 
      inequality in \eqref{eq:problems} holds;
    \item\label{prob:cautious} (One-shot cautious suboptimization)
      given a set $\calS\subseteq \mathbb{R}^n$, find $z\in\calS$ for
      which Problem~\ref{prob:cautious-opt}\ref{prob:verif} yields
      the minimal value of $\delta$;
    \item\label{prob:setwise} (Set-wise verification of suboptimality)
      given a set $\calS\subseteq \mathbb{R}^n$, find the minimal
      $\delta\in\mathbb{R}$ as a function of the data $(Y,\Phi)$, for
      which either inequality 
      in~\eqref{eq:problems} holds for all $z\in\calS$.
    \end{enumerate}
  }
\end{problem}

Given a set of basis functions and a noise model, being able to
resolve any of these problems can be viewed as a property of the
measurements. This is essentially the viewpoint taken with respect to
data-driven control within the data informativity framework,
e.g.,~\cite{HJVW-JE-HLT-MKC:20,HJVW-MKC-JE-HLT:22}: given
$\delta\in\mathbb{R}$ and $\calS\subseteq\mathbb{R}^n$, one could say
that the data $(Y,\Phi)$ is \emph{informative for
  $\delta$-suboptimization on $\calS$} if there exists $z\in\calS$
such that $\norm{\hat{\phi}(z)}_2\leq \delta$.}

We rely on the data informativity interpretation throughout our
technical discussion in Section~\ref{sec:cautious-subopt} to solve 
Problem~\ref{prob:cautious-opt}. Then, in
Section~\ref{sec:applications}, we illustrate how to apply our results
to solve various problems in data-driven analysis and control.
First, we consider data-based contraction analysis for a
class of nonlinear systems parameterized in terms of given basis
functions. This essentially requires us to characterize one-sided
Lipschitz constants for the unknown dynamics in terms of measured
data. In terms of data-driven control, we consider the problem of
data-based suboptimal regulation of unknown linear systems.  Finally, 
in Section~\ref{sec:online}
we analyze online suboptimization, where we
iteratively optimize and collect new measurements to find the
optimizer.

\vspace{-0.8em}\section{Cautious optimization}\label{sec:cautious-subopt}

\vspace{-0.3em}{\color{black} This section describes our solutions to Problem~\ref{prob:cautious-opt}. 
Recall that, based on the data we can perform set-valued regression to obtain the set 
$\Theta$. In fact, we cannot distinguish $\hat{\phi}$ from any of the
functions $\phi^\theta$ as long as $\theta \in \Theta$.} This means
that we can only guarantee that $\norm{\hat{\phi}(z)}_2 \leq \delta$
holds if $\sup_{\theta\in\Theta} \norm{\phi^\theta(z)}_2 \leq
\delta$. Similarly, given a vector $c\in\mathbb{R}^m$, we can only
conclude that $c^\top \hat{\phi}(z) \leq \delta$ if
$\sup_{\theta\in\Theta} c^\top\phi^\theta(z) \leq \delta$. This
motivates the following definitions:
\vspace{-0.4em}\begin{equation}\label{eq:def g gc}
  \vspace{-0.6em}g(z;\Theta):=
  \sup_{\theta\in\Theta} \norm{\phi^\theta(z)}_2,\quad g_c(z;\Theta)
  := \sup_{\theta\in\Theta} c^\top\phi^\theta(z).
\end{equation}
These functions correspond to the elementwise worst-case realization
of the unknown parameter~$\hat{\theta}$. Therefore, we refer to
$g(z;\Theta)$ as the \textit{worst-case norm bound} and
$g_c(z;\Theta)$ as the \textit{worst-case linear bound}.

To draw conclusions regarding the true function $\hat{\phi}$, we
investigate the functions~$g(z;\Theta)$ and~$g_c(z;\Theta)$. In fact,
problems (i)-(iii) in Problem~\ref{prob:cautious-opt} can be recast as
follows.  Resolving Problem~\ref{prob:cautious-opt}\ref{prob:verif}
(in Section~\ref{sec:verification}) boils down to finding
function values of $g(z;\Theta)$ or $ g_c(z;\Theta)$. Resolving
Problem~\ref{prob:cautious-opt}\ref{prob:cautious}
(in Section~\ref{sec:subopt}) is equivalent to finding
\vspace{-0.4em}\begin{equation}\label{eq:min min}
\vspace{-0.6em}  \min_{z\in \calS} g(z;\Theta) \quad
  \textrm{ or }\quad \min_{z\in \calS} g_c(z;\Theta).
\end{equation}
Finally, resolving Problem~\ref{prob:cautious-opt}\ref{prob:setwise}
(in Section~\ref{sec:set-wise}) amounts to finding the minimal
$\delta$ for which
\vspace{-0.4em}\begin{equation}\label{eq:max max}
\vspace{-0.6em}  \max_{z\in \calS} g(z;\Theta)\leq \delta \quad
  \textrm{ or }\quad \max_{z\in \calS} g_c(z;\Theta)\leq \delta.
\end{equation}
{\color{black} Given the basis functions, each of the previous problems
  is completely specified in terms of the set $\Theta$. In turn, this
  set is determined by the matrix $N$, containing information on the
  noise model and the measured data.}  

\vspace{-0.4em}\subsection{Pointwise verification of
  suboptimality}\label{sec:verification}
\vspace{-0.4em}In order to resolve Problem~\ref{prob:cautious-opt}\ref{prob:verif},
we are interested in obtaining values of the functions in
\eqref{eq:def g gc}. As a first observation, note that if the set
$\Theta$ is compact (cf. Remark~\ref{rem:exciting}), then the supremum
is attained, and we can replace it with a maximum. In this case, the
worst-case norm bound and worst-case linear bound are finite-valued
functions.

We rely on the following characterization to obtain function values
of~$g(z;\Theta)$.

\begin{lemma}\longthmtitle{Data-based conditions for bounds on the
    norm}\label{lem:lmi for subopt}
  Given data $(Y,\Phi)$, let $\Theta=\calZ(N)$, with $N$ as in
  \eqref{eq:defN}.  Suppose that $N$ has at least one positive
  eigenvalue.  Let $\delta\geq 0$ and $z\in \mathbb{R}^n$. Then,
  \vspace{-0.3em}\begin{equation}\label{eq:probs norm}
    \vspace{-0.3em}\norm{\phi^\theta(z)}_2 \leq
    \delta \quad \textrm{ for all } \quad
    \theta\in\Theta
  \end{equation}
  if and only if there exists $\alpha\geq 0$ such that
  \vspace{-0.3em}\begin{equation}\label{eq:LMI cond norm}
    \vspace{-0.3em}{\scriptsize\begin{pmat}[{.|}]
      \delta^2 I_m &0 &0 \cr
      0& 0 &b(z) \cr\-
      0& b(z)^\top & 1 \cr
    \end{pmat}} - \alpha  
    \begin{pmat}[{|}]
      N & 0 \cr\- 0 & 0 \cr
    \end{pmat} \geq0 .
  \end{equation}
\end{lemma}
\begin{pf} 
  Note that
  $ \norm{\phi^\theta(z)}^2_2 = b(z)^\top \theta\theta^\top b(z)$.
  This allows us to write \eqref{eq:probs norm} equivalently as
  \vspace{-0.4em}\[\vspace{-0.4em}
    b(z)^\top \theta\theta^\top b(z) \leq \delta^2 \iff \theta^\top
    b(z)b(z)^\top \theta\leq \delta^2 I_m,
  \]
  which can also be expressed as
  \vspace{-0.6em}\begin{equation}\label{eq:norm as qmi}
    \vspace{-0.3em}\scriptsize\begin{bmatrix} I_m \\
      \theta \end{bmatrix}^\top \begin{bmatrix} \delta^2 I_m & 0 \\ 0
      & -b(z)b(z)^\top \end{bmatrix} \begin{bmatrix} I_m \\
      \theta \end{bmatrix}\geq 0.
  \end{equation}
  Thus \eqref{eq:probs norm} holds if and only if \eqref{eq:norm as
    qmi} holds for all $\theta\in\Theta$. Given that
  $\Theta=\calZ(N)$, where $N$ has at least one positive eigenvalue,
  the statement follows from applying~\cite[Thm
  4.7]{HJVW-MKC-JE-HLT:22} to \eqref{eq:norm as qmi} and using a Schur
  complement. \QED
\end{pf}
\vspace{-1em}The necessary and sufficient conditions of Lemma~\ref{lem:lmi for
  subopt} allow us to provide an alternative expression
for~$g(z;\Theta)$ that allows us to compute its value efficiently.

\begin{corollary}\longthmtitle{Values of the worst-case norm bound}
	\label{cor:sup of norm using lmi}
  Given data $(Y,\Phi)$, let $\Theta=\calZ(N)$, with $N$ as in
  \eqref{eq:defN}. Suppose that $N$ has at least one positive
  eigenvalue. Then, for $z\in \mathbb{R}^n$, we have
  \vspace{-0.4em}\begin{equation}\label{eq:value of g}
   \vspace{-0.4em} g(z;\Theta) = \min\{\delta\geq
    0 \mid \exists \alpha\geq 0 \textrm{ s.t. } \eqref{eq:LMI cond
      norm} \textrm{ holds} \}.
  \end{equation}
\end{corollary}

{\color{black} Technically the statement~\eqref{eq:LMI cond norm} is
  not an LMI, since it depends quadratically on $\delta$. However,
  function values of $g(z;\Theta)$ can be found efficiently via
  Corollary~\ref{cor:sup of norm using lmi} by minimizing $\delta^2$
  over \eqref{eq:LMI cond norm} instead. Similar situations will arise
  without further mention in the remainder of the paper. }

\begin{remark}\longthmtitle{Weakening assumptions on data lead to
    conservative upper bound}\label{rem:weakened}
  {\rm If $N$ does \textit{not} have positive eigenvalues, i.e.,
    $N\leq 0$, then the `if'-side of Lemma~\ref{lem:lmi for subopt}
    remains true. Therefore, we can conclude that \eqref{eq:value of
      g} in Corollary~\ref{cor:sup of norm using lmi} holds with
    `$\leq$' replacing `$=$'. This means that we obtain an upper bound
    for $g(z;\Theta)$, leading to a conservative upper bound of
    $\hat{\phi}(z)$. \oprocend}
\end{remark}

Moreover, we can use Corollary~\ref{cor:sup of norm using lmi} to find
bounds independent of the choice of~$z$.

\begin{remark}\longthmtitle{Relaxations of the conditions for bounds
    on the norm}\label{rem:overestimating}
  {\rm If $M$ is such that $b(z)b(z)^\top \leq M$, then
    \[
      {\scriptsize\begin{bmatrix}
        \delta^2 I_m
        & 0
        \\
        0 & -b(z)b(z)^\top
      \end{bmatrix}}
      \geq
      {\scriptsize\begin{bmatrix}
        \delta^2 I_m & 0
        \\
        0 & -M
      \end{bmatrix}}.
    \]
    This implies that 
\vspace{-0.5em}    \begin{equation}\label{eq:cond in norm of b}
\vspace{-0.5em}      {\scriptsize\begin{bmatrix}
        \delta^2 I_m & 0
        \\
        0 & -M
      \end{bmatrix}}-\alpha
      N \geq 0
    \end{equation}
    implies \eqref{eq:LMI cond norm}. In turn, this means that for any
    $M$ such that $b(z)b(z)^\top \leq M$ for all $z$, we have
    \[
      g(z;\Theta) \leq \min\{\delta\geq 0 \mid \exists \alpha\geq 0
      \textrm{ s.t. }  \eqref{eq:cond in norm of b} \textrm{ holds}
      \},
    \]
    which yields a bound independent from~$z$.
    \oprocend } 
\end{remark}  

We rely on the following characterization to obtain function values
of~$g_c(z;\Theta)$.

\begin{lemma}\longthmtitle{Data-based conditions for scalar upper
    bounds}\label{lem:lmi
    for linear} Given data $(Y,\Phi)$, let $\Theta=\calZ(N)$, with $N$
  as in \eqref{eq:defN}. Let $\delta\in\mathbb{R}$,
  $z\in \mathbb{R}^n$, and $c\in\mathbb{R}^m$. Suppose
  $c^\top(N| N_{22})c>0$ and let
  \vspace{-0.4em}\begin{equation}\label{eq:def Nc} \vspace{-0.4em}N_c
    :=
    {\scriptsize\begin{bmatrix}
      c^\top
      N_{11}c & c^\top N_{12}
      \\ N_{21}c &
                   N_{22}
    \end{bmatrix}}.
  \end{equation}
  Then it holds that
  \vspace{-0.5em}\begin{equation}\label{eq:probs value}
    \vspace{-0.5em}c^\top\phi^\theta(z) \leq \delta \quad \textrm{ for all } \quad
    \theta\in\Theta  
  \end{equation}
  if and only if there exists $\alpha\geq 0$ such that
 \vspace{-0.5em} \begin{equation}\label{eq:LMI cond value}
   \vspace{-0.7em} {\scriptsize\begin{bmatrix} 2\delta &-b(z)^\top\\-b(z)&0\end{bmatrix}}- \alpha
    N_c \geq 0. 
  \end{equation}
\end{lemma}
\vspace{-1.3em}\begin{pf}
  We can rewrite $c^\top\phi^\theta(z)$ as
  $ c^\top\theta^\top b(z) + b(z)^\top \theta c \leq 2 \delta$.  Thus,
  \eqref{eq:probs value} holds if and only if
  \vspace{-0.5em}\[\vspace{-0.5em}
    \gamma^\top b(z) +b(z)^\top \gamma \leq 2\delta \quad \textrm{ for
      all } \quad \gamma \in \calZ(N)c,
  \]
  where $\calZ(N)c := \{Zc \mid Z\in\calZ(N)\}$. Since $c\neq 0$,
  by~\cite[Thm. 3.4]{HJVW-MKC-JE-HLT:22}, we have
  $\calZ(N)c = \calZ(N_c)$. Therefore, \eqref{eq:probs value} is
  equivalent to: \vspace{-0.4em}
  \begin{equation}\label{eq:probs value
      2} \vspace{-0.4em}
    {\scriptsize
      \begin{bmatrix}
        1 \\
        \gamma
      \end{bmatrix}^\top
      \begin{bmatrix}
      2\delta
      &-b(z)^\top
      \\
      -b(z)&0
    \end{bmatrix}
    \begin{bmatrix}
      1
      \\
      \gamma
    \end{bmatrix}}
    \geq 0, \textrm{ for all }
    \gamma\in\calZ(N_c).
  \end{equation}
  Given that $N_{22}\leq 0$, we have that $c^\top(N| N_{22})c>0$ if
  and only if $N_c$ has a positive eigenvalue. Then, the statement
  follows from applying~\cite[Thm 4.7]{HJVW-MKC-JE-HLT:22}
  to~\eqref{eq:probs value 2}. \QED
\end{pf}
\vspace{-0.5em}%
We leverage the characterization in Lemma~\ref{lem:lmi for linear} to
provide a closed-form expression of $ g_c(z;\Theta)$ in terms of the
data.

\begin{theorem}\longthmtitle{Explicit bounds in the scalar
    case}\label{thm:explicit delta} Given data $(Y,\Phi)$, let
  $\Theta=\calZ(N)$, with $N$ as 
  in \eqref{eq:defN}. Let $z\in \mathbb{R}^n$ and $c\in\mathbb{R}^m$. 
  If $b(z)\in \im\Phi$ or $c=0$, then
\vspace{-0.3em}  \[\vspace{-0.6em}	\begin{split}
g_c(z;\Theta)=& -c^\top N_{12}N_{22}^\dagger b(z)
	\\
	&\quad +
	\sqrt{c^\top(N|N_{22})cb(z)^\top (-N_{22}^\dagger) b(z)},  	
\end{split}\]
  and, if $c\neq 0$ and $b(z)\not\in \im\Phi$, then
  $g_c(z;\Theta)= \infty$.
\end{theorem}
\vspace{-0.8em}\begin{pf}
  Note that, if either $c=0$ or $b(z)=0$, then we have
  $g_c(z;\Theta) = 0$, thus the result follows. Assume
  then that $c\neq 0$, $b(z)\neq 0$, and $b(z)\in \im\Phi$.  If
  $c^\top(N| N_{22})c= 0$, then $N_c\leq 0$. As such,
  $\gamma\in\calZ(N_c)$ if and only if
\vspace{-0.4em}  \[
\vspace{-0.6em}    N_c
    {\scriptsize \begin{bmatrix}
      1
      \\
      \gamma
    \end{bmatrix}}
    = 0 \iff \gamma \in -N_{22}^\dagger N_{21}c+\ker N_{22}.
  \]
  Using \eqref{eq:probs value 2}, we see that
  $c^\top \phi^\theta\leq \delta$ for all $\theta\in\Theta$ if and
  only if $\gamma^\top b(z)\leq \delta$ for all
  $\gamma \in \calZ(N_c)$. Combining the previous, we see that this
  holds if and only if $\delta\geq -c^\top N_{12}N_{22}^\dagger b(z)$
  (where we have used that $b(z)\in \im\Phi$). Minimizing $\delta$
  reveals that the statement holds for this case.

  Consider the case when $c^\top(N| N_{22})c>0$. By Lemma~\ref{lem:lmi
    for linear}, \eqref{eq:probs value} holds if and only if
  \eqref{eq:LMI cond value} holds for some $\alpha \ge 0$. Since
  $b(z)\neq 0$, it must be that $\alpha >0$. Using the Schur
  complement on \eqref{eq:LMI cond value}, we obtain
\vspace{-0.7em}  \begin{align*}
  0& \leq\! 2\delta\!-\!\alpha c^{\!\top} N_{11}c +\alpha^{\!-1} (b(z)^{\!\top} +\alpha
       c^{\!\top} N_{12})N_{22}^\dagger (b(z) \!\! +\alpha N_{21}c)
    \\
     & = 2\delta + c^\top N_{12}N_{22}^\dagger b(z) +b(z)^\top
       N_{22}^\dagger N_{21}c
    \\ 
     &\quad - \alpha c^\top (N|N_{22})c +\alpha^{-1} b(z)^\top N_{22}^\dagger b(z). 
  \end{align*}
  Clearly, there exists $\alpha$ for which this holds if and only if
  it holds for the value of $\alpha$ that maximizes the expression on
  the right. Note that $\ker N_{22} = \ker\Phi^\top$. Thus, since
  $b(z)\in \im \Phi$ and $b(z)\neq0$, we see that
  $b(z)^\top N_{22}^\dagger b(z)<0$.  The expression is then maximized
  over $\alpha\geq 0$ with
\vspace{-0.3em}  \[
\vspace{-0.3em}    \alpha = {\scriptsize\sqrt{\frac{-b(z)^\top N_{22}^\dagger b(z)}{c^\top
        (N|N_{22})c }}}.
  \]
  This yields that \eqref{eq:LMI cond value} holds if and only if
\vspace{-0.8em} \begin{align*}
  &2\delta + c^\top N_{12}N_{22}^\dagger b(z) +b(z)^\top
      N_{22}^\dagger N_{21} c
    \\
    &- 2 \sqrt{ c^\top(N|N_{22})c b(z)^\top
      (-N_{22}^\dagger) b(z)} \geq 0.
  \end{align*}
  To obtain $g_c (z;\Theta)$, we minimize this over $\delta$, which
  yields the expression in the statement.

  Finally, consider the case when $c\neq 0$ and
  $b(z)\not\in \im \Phi$. It is straightforward to check that if
  $\theta\in\Theta$, then $\theta+vw^\top \in \Theta$, for all
  $w\in\mathbb{R}^m$ and $v\in\ker \Phi^\top \subseteq
  \mathbb{R}^k$. Since $b(z)\not\in \im\Phi$ and $\im \Phi$ is
  orthogonal to $\ker \Phi^\top$, we can take $v$ such that
  $v^\top b(z)>0$. Now, for arbitrary $\beta\in\mathbb{R}$, take
  $w=\beta c$. Then,
  \vspace{-0.3em}  \[\vspace{-0.3em}  
    \sup_{\theta\in\Theta}\phi^\theta(z) \geq c^\top(\theta
    +vw^\top)^\top b(z)= c^\top\theta^\top b(z) +\beta c^\top cv^\top
    b(z).
  \]
  Since $c^\top c>0$ and $v^\top b(z)>0$, and this is valid for any
  $\beta\in\mathbb{R}$, the last part of the statement follows.  \QED
\end{pf}

\begin{remark}[The special case of $N_{22}<0$]\label{rem:N22 less than
    0}%
  {\rm Note that when $\Theta$ is compact (or equivalently, $\Phi$
    has full row rank, cf. Remark~\ref{rem:exciting}), it is immediate
    that $b(z)\in\im\Phi$ for any $z\in\mathbb{R}^n$, and hence
    $g_c(z;\Theta)$ is always finite valued. \oprocend}
\end{remark}

The following illustrates the introduced concepts in the familiar
setting of linear regression.

\begin{example}[Linear regression]\label{ex:linear}
  {\rm Suppose that $\hat{\phi}:\mathbb{R}^2\rightarrow\mathbb{R}^2$
    and $k=3$. We consider basis functions $\phi_1(z) =1$ and
    $\phi_2(z)=\begin{pmatrix} 1 & 0 \end{pmatrix} z$,
    $\phi_3(z) = \begin{pmatrix} 0 & 1 \end{pmatrix}z$. In other
    words, we take $\theta\in \mathbb{R}^{3\times 2}$, and let 
    $\phi^\theta (z) = \theta^\top
      \begin{bmatrix}
        1 &z^\top
      \end{bmatrix}^{\!\top}$. 
    We collect three measurements,
   \vspace{-0.4em}  \[ \vspace{-0.4em}
      \begin{pmatrix}
        1
        \\
        1
      \end{pmatrix}
      =
      \hat{\phi}(0) +w_1, 
      \begin{pmatrix}
        0
        \\
        1
      \end{pmatrix}
      = \hat{\phi}(e_1) +w_2, 
      \begin{pmatrix}
        1
        \\
        0
      \end{pmatrix}
      = \hat{\phi}(e_2) +w_3,
    \]
    where we assume that $WW^\top\leq I_2$. 
    This leads to $\Theta=\calZ(N)$ with $N$ such that
    \vspace{-0.6em}\[
      N_{11} = {\scriptsize\begin{bmatrix} 1 & 0 \\ 0 & 1 \end{bmatrix}}, 
      N_{12} = {\scriptsize\begin{bmatrix}2 &0& 1 \\ 2& 1 &0\end{bmatrix}}, 
      N_{22} =
      {\scriptsize\begin{bmatrix}
        -3 & -1& -1
        \\
        -1& -1 & 0
        \\
        -1& 0&-1
      \end{bmatrix}} .\vspace{-0.5em}
    \]
    Now, the least-squares estimate is equal to
    \vspace{-0.5em}\[
      \theta^{\lse} = 
      {\scriptsize\begin{bmatrix} 1& 1\\ -1&0 \\0 &-1 \end{bmatrix}},
      \quad      	
      	\phi^{\lse}(z;\Theta)=
      {\scriptsize\begin{pmatrix}
        1 -z_1
        \\
        1 -z_2
      \end{pmatrix}}.\vspace{-0.5em}
    \]
    In fact, using Theorem~\ref{thm:explicit delta} we can conclude
    that
  \vspace{-0.4em}  \[ \vspace{-0.4em}
      g_c(z;\Theta) = c^\top\!\!
      {\scriptsize\begin{pmatrix}
        1-z_1 \\
        1-z_2
      \end{pmatrix}}
      + \sqrt{c^\top
        c}\sqrt{(1-z_1-z_2)^2+z_1^2+z_2^2}.
    \]
    As such, on the basis of the data we can guarantee that, for
    instance,
    $c^\top\hat{\phi}(0)\leq (1 \hspace{1em} 1)c+ \norm{c}_2$. \oprocend }
\end{example}

\begin{figure}[t]
	\begin{center}
		\includegraphics[trim=0cm 4cm 0cm 4cm,clip,width=8.6cm]{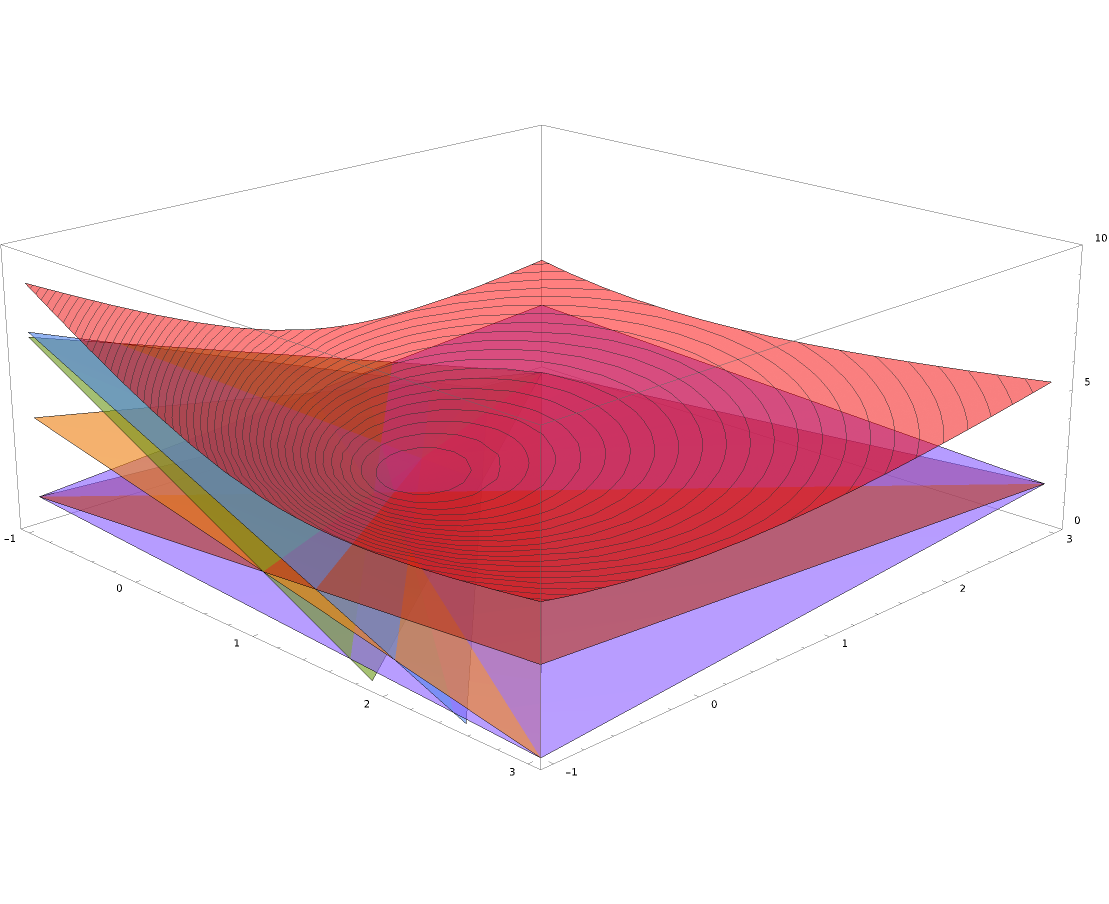}
		\includegraphics[width=4.3cm]{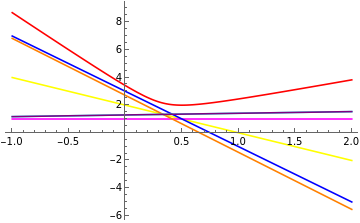}
		\includegraphics[width=4.3cm]{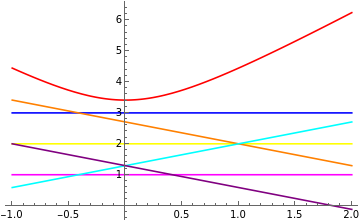}
	\end{center}
	\caption{A plot depicting the situation of
          Example~\ref{ex:linear}. For
          $c=\left(\begin{smallmatrix} 1
              \\1 \end{smallmatrix}\right)$, we plot a number function
          $c^\top\phi^\theta$, with $\theta$ consistent with the
          data. Any of these functions could be the unknown function.
          The least squares estimate is shown in yellow. If we are
          interested in bounds on $\hat{\phi}$, we can employ
          $g_c(z;\Theta)$ (shown in red).  In particular note that
          this is not a linear map. Below the plot, we show the slices
          $z_1=z_2$ and $z_1=-z_2$ respectively.}
	\label{fig:plot example}
\end{figure}

\subsection{Uncertainty on function values}\label{sec:uncertainty}
The results in Section~\ref{sec:verification} allow us to find
bounds for the unknown function $\hat{\phi}$, but do not consider
how much these bounds deviate from its true value. For this, note that 
the result of Theorem~\ref{thm:explicit delta} can be interpreted 
as quantifying the
distance between the true function value and the least-squares
estimate $\phi^{\lse}$ (defined in Remark~\ref{rem:lse}) in terms of 
the basis functions and the data. This observation leads us to define the
\textit{uncertainty} corresponding to~\eqref{eq:def g gc} at $z$ by
 \vspace{-0.5em}\begin{subequations} \vspace{-0.2em}
  \begin{align}
    U(z;\Theta)
    &:= \sup_{\theta\in\Theta}
      \norm{\phi^\theta(z)-\phi^{\lse}(z;\Theta)}_2,
    \\ 
    U_c(z;\Theta)
    &:=  \sup_{\theta\in\Theta}
      c^\top\big(\phi^\theta(z)-\phi^{\lse}(z;\Theta)\big). 
\end{align}
\end{subequations}  
Thus, if the uncertainty $U(z;\Theta)$ at $z$ is small, then
$\phi^{\lse}(z;\Theta)$ is quantifiably a good estimate of $\hat{\phi}(z)$.

The uncertainty also allows us describe a useful variant to
Problem~\ref{prob:cautious-opt}\ref{prob:cautious}. In order to
balance the demands of a low upper bound on the value of the true
function with an associated low uncertainty, it is reasonable to
consider the following generalization of the cautious suboptimization
problem in~\eqref{eq:min min}: for $\lambda\geq 0$, consider
 \vspace{-0.4em}\begin{subequations} 
  \begin{align}
    &\min_{z\in\calS} g(z;\Theta) +\lambda U(z;\Theta), \label{eq:weighted cs a}
    \\ 
    &\min_{z\in \calS} g_c(z;\Theta)+\lambda
      U_c(z;\Theta). \label{eq:weighted cs b} 
  \end{align}
\vspace{-1.8em}\end{subequations}

The following result provides a closed-form expression for the
objective function in~\eqref{eq:weighted cs a}.

\begin{lemma}\longthmtitle{Explicit expression of uncertainty in norms}
  Given data $(Y,\Phi)$, let $\Theta=\calZ(N)$, with $N$ as
  in~\eqref{eq:defN}. For $z\in \mathbb{R}^n$, if $b(z)\in \im\Phi$,
  then
   \vspace{-0.7em}\begin{equation}\label{eq:expl U}
    \vspace{-0.4em} U(z;\Theta) =
    \sqrt{\lambda_{\operatorname{max}}(N|N_{22})
      b(z)^\top(-N_{22}^\dagger)b(z)}.
  \end{equation}
  If $b(z)\not\in\im \Phi$, then $U(z;\Theta)= \infty$.
\end{lemma}
 \vspace{-0.4em}\begin{pf}
  First suppose that $N$ has at least one positive eigenvalue. Similar
  to the proof of Lemma~\ref{lem:lmi for subopt}, we can conclude that
  that
  $\norm{ \phi^\theta(z) - \phi^{\lse}(z;\Theta)}_2\leq \varepsilon$
  for $\varepsilon\geq 0$ and for all $\theta\in\Theta$ if and only if
  there exists $\alpha\geq0$ such that
  \begin{equation}\label{eq:LMI lse thing}
    {\scriptsize\begin{bmatrix}
      \varepsilon^2 I & 0\\ 0 & -b(z)b(z)^\top \end{bmatrix}}
    -\alpha {\scriptsize\begin{bmatrix} N|N_{22} & 0\\ 0 & N_{22}\end{bmatrix}}\geq
    0.
  \end{equation}
  Thus, in line with Corollary~\ref{cor:sup of norm using lmi}, we obtain
  \begin{equation}\label{eq:U as lmi}
    U(z;\Theta) = \min\{\varepsilon\geq 0 \mid \exists \alpha\geq 0
    \textrm{ s.t. } \eqref{eq:LMI lse thing} \}.
  \end{equation}
  If $b(z)\in \im \Phi$ then we can minimize $\varepsilon$ by taking
  $\alpha = b(z)^\top(-N_{22}^\dagger)b(z)$, leading to \eqref{eq:expl
    U}. If $N$ does not have at least one positive eigenvalue, then
  $N|N_{22}=0$. As in Remark~\ref{rem:weakened}, we can conclude
  \eqref{eq:U as lmi} with an inequality instead. Moreover, it
  immediately follows that $U(z;\Theta)\leq 0$, proving that
  \eqref{eq:expl U} holds.  The second part of the statement can be
  concluded in a similar fashion as the last part of the proof of
  Theorem~\ref{thm:explicit delta}. \QED
\end{pf}
\vspace{-1em}In addition to expressing the uncertainty,
this result also gives rise to a useful explicit bound on the function
$g(z;\Theta)$.

\begin{corollary}\longthmtitle{Explicit upper bound of the
    worst-case norm bound}\label{cor:upper g}
  Given data $(Y,\Phi)$, let $\Theta=\calZ(N)$, with $N$ as in
  \eqref{eq:defN}. For $z\in \mathbb{R}^n$, if $b(z)\in \im \Phi$,
  then
\vspace{-0.4em}  \begin{equation}\label{eq:expl g}\vspace{-0.8em}
  	\begin{split}
    g(z;\Theta) &\leq
    \norm{\phi^{\lse}(z;\Theta)}_2
    \\
    &\quad+
    \sqrt{\lambda_{\operatorname{max}}(N|N_{22})
      b(z)^\top(-N_{22}^\dagger)b(z)}.
  \end{split}\end{equation}
\end{corollary}
The following result provides a closed-form expression for the
objective function in~\eqref{eq:weighted cs b}.

\begin{lemma}\longthmtitle{Explicit expression of scalar
    uncertainty}\label{lem:expl uncertainty}
  Given data $(Y,\Phi)$, let $\Theta=\calZ(N)$, with $N$ as in
  \eqref{eq:defN}. Let $z\in \mathbb{R}^n$ and $c\in\mathbb{R}^m$.  If
  $b(z)\in \im\Phi$ or $c=0$, then
  \vspace{-0.4em}\[
    U_c(z;\Theta)= \sqrt{(c^\top(N|N_{22})c)(b(z)^\top
      (-N_{22}^\dagger) b(z))}.
  \]
  If $b(z)\not\in \im\Phi$ and $c\neq0$, then $U_c(z;\Theta)=\infty$. 
  Regardless, if $\lambda\geq 0$ and we define 
  \vspace{-0.4em}\[
  \vspace{-0.3em} N_\lambda := {\scriptsize\begin{bmatrix} N_{11} & N_{12} \\ N_{21} &
      N_{22} \end{bmatrix}}+ {\scriptsize\begin{bmatrix}
      \lambda(2+\lambda)(N|N_{22}) & 0\\ 0&0 \end{bmatrix}},
  \]
  then $   g_c(z;\Theta)+\lambda U_c(z;\Theta) = g_c(z;\calZ(N_\lambda))$.
\end{lemma}
\vspace{-0.7em}\begin{pf}
  The first part follows immediately from Theorem~\ref{thm:explicit
    delta}. Now note that
  \vspace{-0.3em}\[
   \vspace{-0.3em} N_\lambda = {\scriptsize\begin{bmatrix}
      (1+\lambda)^2(N|N_{22})+N_{12}N_{22}^{-1}N_{21} & N_{12} \\
      N_{21} & N_{22} \end{bmatrix}}.
  \]
  Therefore $N_\lambda|N_{22} = (1+\lambda)^2(N|N_{22})$. Thus, the
  result follows from the application of Theorem~\ref{thm:explicit
    delta}.\QED
\end{pf}
\vspace{-1em}Lemma~\ref{lem:expl uncertainty} means that, even though
problem~\eqref{eq:weighted cs b} is more general than the
corresponding cautious suboptimization problem in~\eqref{eq:min min},
both problems can be resolved in the same fashion.

\subsection{Convexity and suboptimization}\label{sec:subopt}

{\color{black} Solving the one-shot cautious suboptimization problem,
  cf.  Problem~\ref{prob:cautious-opt}\ref{prob:cautious}, in an
  efficient manner, amounts to minimizing the functions
  $g(\cdot;\Theta)$ or $g_c(\cdot ;\Theta)$, respectively. For this,
  any standard optimization method can be used, which we illustrate
  below using simple gradient descent.  We discuss how to determine
  values of the Jacobians of $g(\cdot;\Theta)$ or $g_c(\cdot ;\Theta)$
  in Appendix~\ref{sec:appendix}, where we give a simple
  characterization in terms of properties of the basis and data, cf.
  Corollary~\ref{cor:gradient}. Under suitable regularity conditions,
  these Jacobians can be used in a gradient descent scheme, which
  converges to a local minimum.

  We can also resolve
  Problem~\ref{prob:cautious-opt}\ref{prob:cautious} \textit{globally}
  in this manner if $g(\cdot;\Theta)$ or $g_c(\cdot ;\Theta)$ are
  convex, respectively. To provide efficient tests for convexity, we
  will first assume that $\theta c$ is elementwise nonnegative for all 
  $\theta\in\Theta$. A test for this property can be found in the appendix,
  in Lemma~\ref{lem:theta positive}. Using this, we can 
  identify the following conditions which ensure the
	convexity of $\norm{\hat{\phi}}_2$ or
  $c^\top\hat{\phi}$ and their upper bounds.} 

\vspace{-0.2em}\begin{proposition}\longthmtitle{Convexity of the true
    function}\label{cor:convexity}
  Given data $(Y,\Phi)$, let $\Theta=\calZ(N)$, with $N$ is as in
  \eqref{eq:defN}. Assume the basis functions $\phi_i$ are convex.
  Let $c\neq 0$ with $\theta c\in\mathbb{R}^k_{\geq 0} $ for all
  $\theta\in\Theta$. Then,
  \begin{itemize}
  \item $c^\top\phi^\theta$ is convex for all $\theta\in\Theta$ and
    $g_c(\cdot;\Theta)$ is a finite-valued convex function;
  \item If, in addition, the functions $\phi_i$ are strictly convex
    and $0\not\in \calZ(N_c)$, then $c^\top\phi^\theta$ is strictly
    convex for all $\theta\in\Theta$ and $g_c(\cdot;\Theta)$ is
    strictly convex.
  \end{itemize}
  Moreover, $\norm{\phi^\theta(\cdot)}_2$ is convex if
  $e_i^\top\phi^\theta$ is convex and nonnegative for all
  $i=1,\ldots m$.
\end{proposition}
\vspace{-1.5em}\begin{pf}
  Lemma~\ref{lem:theta positive} shows that, if
  $\theta c\in\mathbb{R}^k_{\geq 0} $ for all $\theta\in\Theta$, then
  $\Phi$ has full row rank or, equivalently, $N_{22}<0$. This implies,
  cf. Remark~\ref{rem:exciting}, that $\calZ(N)$ is compact and
  therefore $g_c(\cdot;\Theta)$ is finite-valued. The first result now
  follows from the facts that nonnegative combinations of convex
  functions are convex and that the maximum of any number of convex
  functions is also convex. The last part follows from standard
  composition rules for convexity, see e.g.,~\cite[Example
  3.14]{SB-LV:09}.\QED
\end{pf}
\vspace{-1.4em}

{\color{black} The results of Lemma~\ref{lem:theta positive} and
  Proposition~\ref{cor:convexity} taken together mean that, if the basis
  functions are convex, we can test for convexity on the basis of
  data. To guarantee strict convexity, Proposition~\ref{cor:convexity}
  requires that $0\not\in\calZ(N_c)$, which is easy to check in terms
  of data. }

Looking back at Example~\ref{ex:linear}, we can make two interesting
observations.  First, the conditions of Proposition~\ref{cor:convexity}
do \textit{not} hold in the example, since for instance
$\theta^{\lse}c\in\mathbb{R}^3_{\geq 0}$ only if $c=0$.  Yet
$c^\top\phi^\theta$ is convex for all $\theta\in\Theta$. This is
because the basis functions $\phi_i$ are linear, and therefore the
coefficients $\theta_i$ are not required to be nonnegative. Second,
$g_c(\cdot;\Theta)$ in Example~\ref{ex:linear} is strictly convex for
$c\neq 0$ even though the basis functions are not. Recall that we are
interested in the optimization problem~\eqref{eq:min min}, and
thus not necessarily in properties of the true function
$\hat{\phi}$, but of the worst-case linear bound
$g_c(\cdot;\Theta)$. These observations motivate our ensuing
discussion to provide conditions for convexity of the upper
bound instead. Under the assumptions of
Theorem~\ref{thm:explicit delta}:
\vspace{-0.6em}\[\vspace{-0.6em}
g_c(z;\Theta) = c^\top\phi^{\lse}(z;\Theta) + U_c(z;\Theta).
\]
Thus, if (i)
$c^\top \phi^{\lse}= -c^\top N_{12}N_{22}^\dagger
b(\cdot)$ is convex and (ii) $U_c$ is convex, then so is
$g_c$. Moreover, if in addition either is strictly
convex, then so is $g_c$. Condition (i) can be checked
directly if all basis functions $\phi_i$ are twice continuously
differentiable by computing the Hessian of $c^\top\phi^{\lse}$.  Here,
we present a simple criterion, also derived from composition rules for
convexity (again, see e.g.,~\cite[Example 3.14]{SB-LV:09}), to test for
condition~(ii).

\vspace{-0.2em}\begin{corollary}\longthmtitle{Convexity of the
    uncertainties}\label{cor:conv uncert}
  Given data $(Y,\Phi)$, let $\Theta=\calZ(N)$ and assume
  $N_{22}<0$. Then both $U(\cdot;\Theta)$ and $U_c(\cdot;\Theta)$ are
  convex if each basis function $\phi_i$ is convex and
  $- N_{22}^{-1}b(z)\geq 0$ for all $z\in\mathbb{R}^n$.
\end{corollary} 

In addition to guaranteeing that $g_c(\cdot;\Theta)$ is convex, this
result can be used to find a convex upper bound when optimizing 
the norm, by bounding $g(\cdot;\Theta)$ using Corollary~\ref{cor:upper g}.

{\color{black} Equipped with the results of this section, one can solve
  Problem~\ref{prob:cautious-opt}\ref{prob:cautious} efficiently using
  the following steps:
\begin{enumerate}
\item Under the assumptions of Theorem~\ref{thm:explicit delta}, we
  can write a closed-form expression for $g_c(\cdot;\Theta)$;
\item We can apply gradient descent using values of the Jacobian found
  with Corollary~\ref{cor:gradient}, which, under suitable regularity
  conditions yields a local minimum;
\item We can test whether this closed form is (strictly) convex using
  e.g., Proposition~\ref{cor:convexity} or Corollary~\ref{cor:conv
    uncert}, and, if so, conclude that the obtained minimum is global.
\end{enumerate}
}

\subsection{Set-wise verification of suboptimality}\label{sec:set-wise}

Here we solve Problem~\ref{prob:cautious-opt}\ref{prob:setwise} and
provide upper (and lower) bounds as in \eqref{eq:problems} which are
guaranteed to hold for all $z\in\calS\subseteq\mathbb{R}^n$. In terms 
of $g$ and $g_c$, this holds only if~\eqref{eq:max max} holds.

{\color{black} To begin, we consider methods on the basis of convexity
  and concavity. First, we see that
   $g_c(z;\Theta)$ is concave if and only if
$g_{-c}(\cdot;\Theta)=-g_c(\cdot;\Theta)$ is convex. Moreover,}
\vspace{-0.4em}\[\vspace{-0.4em}
  \max_{z\in\calS} g_c(z;\Theta) = -\min_{z\in\calS}
  g_{-c}(z;\Theta).
\]
Therefore, we can test for concavity of $g_c$ and apply the
minimization results to $g_{-c}$ of Section~\ref{sec:subopt}.
Beyond the case of concavity, we can also employ convexity in order to
efficiently provide upper bounds over the convex hull of a finite set,
as the following result shows.

\begin{proposition}\longthmtitle{Convexity and maximal
    values}\label{prop:convex for sets}
  Given data $(Y,\Phi)$, let $\Theta=\calZ(N)$.  Let
  $\calF\subseteq\mathbb{R}^n$ be a finite set and let
  $\calS = \conv \calF$. Suppose that $g(z;\Theta)$ is convex. Then
  \eqref{eq:probs norm} for all $z\in\calS$ if and only if $
    g(z;\Theta)\leq \delta \quad \forall z\in\calF$. 
  Similarly, if $g_c(z;\Theta)$ is convex, then 
  \vspace{-0.4em}\[\vspace{-0.4em} \norm{\phi^\theta(z)}_2 \leq
  \delta \quad \textrm{ for all }
  \theta\in\Theta \textrm{ and all } z\in\calS\] 
   if and only if $g_c(z;\Theta)\leq \delta \quad \forall z\in\calF$.
\end{proposition}

The proof of this result follows immediately from the fact that any
convex function attains its maximum over $\calS$ on~$\calF$.  Thus,
convexity (resp., concavity) allows us to find minimal upper bounds
(resp., maximal lower bounds) over a set, leading to a solution of
Problem~\ref{prob:cautious-opt}.\ref{prob:setwise}. 

{\color{black} To solve
  Problem~\ref{prob:cautious-opt}\ref{prob:setwise} in scenarios where
  convexity or concavity cannot be guaranteed, we employ arguments
  based on Lipschitz continuity and coverings.  We are thus interested
  in the question of whether, for all
  $z,z^\star\in\calS\subseteq\mathbb{R}^n$,
\vspace{-0.4em}  \[\vspace{-0.4em}
    \norm{\hat{\phi}(z)-\hat{\phi}(z^\star)}_2 \leq L
    \norm{z-z^\star}_2,
  \]
  holds.  We give a characterization of Lipschitz constants in terms
  of data in Appendix~\ref{sec:appendix}, cf.
  Lemma~\ref{lem:lipschitz}.
  In addition, we say that set $\calG\subset \mathbb{R}^n$ is an
  $\varepsilon$\textit{-covering} of $\calS\subseteq\mathbb{R}^n$ if,
  for every $z\in\calS$, there exists $z^\star\in\calG$ such that
  $\norm{z-z^\star}_2 \leq \varepsilon$. If $\calS$ is a bounded set,
  then for any $\varepsilon>0$, there exists an $\varepsilon$-covering
  with finitely many elements.

  Having access to a Lipschitz constant and an $\varepsilon$-covering
  allows us to bound the unknown function $\hat{\phi}$ as follows.}

\begin{theorem}\longthmtitle{Coverings and bounds}\label{thm:covering}
  Given data $(Y,\Phi)$, let $\Theta=\calZ(N)$. For $\varepsilon>0$,
  let $\calG$ be a finite $\varepsilon$-covering of
  $\calS\subseteq\mathbb{R}^n$. Suppose that $g(z;\Theta)$ is
  Lipschitz with constant $L$ on $\calS$. Then, for all $z\in\calS$,
  \vspace{-0.4em}\[\vspace{-0.4em}
    \min_{z^\star\in\calG} g(z^\star;\Theta) -\varepsilon L \leq
    g(z;\Theta) \leq \max_{z^\star\in\calG} g(z^\star;\Theta)
    +\varepsilon L.
  \]
  If $g_c(z;\Theta)$ is Lipschitz with constant $L$, then
  for all $z\in\calS$,
  \[
    \min_{z^\star\in\calG} g_c(z^\star;\Theta) -\varepsilon L \leq
    g_c(z;\Theta) \leq \max_{z^\star\in\calG} g_c(z^\star;\Theta)
    +\varepsilon L.
  \]
\end{theorem}

\vspace{-0.4em}We omit the proof, which follows directly from combining the
definitions of Lipschitz continuity and $\varepsilon-$coverings.

This result means that we can find guaranteed upper and lower bounds
of either $\norm{\hat{\phi}(z)}$ or $c^\top\hat{\phi}(z)$ over the
bounded set $\calS$ in terms of a finite number of evaluations of
$g(z; \Theta)$ or $g_c(z;\Theta)$, resp. In turn, recall that
Corollary~\ref{cor:sup of norm using lmi} allows us to efficiently
find function values of~$g(z;\Theta)$ on the basis of measurements
and, similarly, Theorem~\ref{thm:explicit delta} allows us to directly
calculate values of~$g_c(z;\Theta)$.

\vspace{-0.7em}\section{Applications to system analysis and control}\label{sec:applications}
In this section, we exploit the proposed solutions to
Problem~\ref{prob:cautious-opt} in two applications: contraction
analysis of unknown nonlinear systems and regulation of an unknown
linear system to a suboptimal point of an unknown cost function.

\subsection{Data-based contraction analysis for nonlinear
  systems}\label{sec:contraction}
Consider the autonomous discrete- or continuous-time system given by
either
\vspace{-0.5em}\begin{equation}\label{eq:nonl system}
\vspace{-0.3em}    z_{k+1} = \hat{\phi}(z_k), \textrm{ or} \quad \dot{z}(t)
    = \hat{\phi}(z(t)),\vspace{-0.3em}
\end{equation}
where $\hat{\phi}: \mathbb{R}^n \rightarrow \mathbb{R}^n$ is an
unknown function. {\color{black} We consider noisy measurements~\eqref{eq:Y-W} as
described in Section~\ref{sec:problem}. Performing set-valued regression 
with this data, we obtain $\hat{\phi}=\phi^\theta$ for some $\theta\in\Theta=\calZ(N)$. In the
continuous-time case, this means that we collect noisy measurements of
the derivative of the state at a finite set of states.  
  In most applications, these derivative measurements need to be
  determined from collected state measurements (see
  e.g.~\cite{JE-JC:23-tac} for a discussion on this assumption).}

The discrete-time system in~\eqref{eq:nonl system} is \textit{strongly
  contracting} (cf.~\cite[Sec. 3.4]{FB:23}) with respect to the
weighted norm $\norm{\cdot}_{2,P^{1/2}}$ if $\hat{\phi}$ admits a
Lipschitz constant $L<1$, that is,
\vspace{-0.3em}\[\vspace{-0.3em}
  \norm{\hat{\phi}(z)-\hat{\phi}(z^\star)}_{2,P^{1/2}} \leq
  L\norm{z-z^\star}_{2,P^{1/2}}.
\]
One can now directly employ Lemma~\ref{lem:lipschitz} in order to
obtain data-based conditions under which this property holds.

To deal with the situation of continuous-time systems, we first
require a number of prerequisites.  Given $P>0$, a function
$f:\mathbb{R}^n\rightarrow\mathbb{R}^n$ is \textit{one-sided Lipschitz
  with respect to} $\norm{\cdot}_{2,P^{1/2}}$ if there exists
$\gamma\in\mathbb{R}$ such that, for all $z,z^\star\in\mathbb{R}^n$,
\vspace{-0.3em}\begin{equation}\label{eq:osl condition}
  (z-z^\star)^\top P
  (f(z)-f(z^\star))\leq
  \gamma\norm{z-z^\star}^2_{2,P^{1/2}}.\vspace{-0.3em}
\end{equation}
Such $\gamma$ is called a \textit{one-sided Lipschitz constant} and
the smallest such $\gamma$ is denoted $\osL(f)$. The autonomous system
$\dot{z}(t)=f(z(t))$ is \textit{strictly contracting} with rate $|b|$
if $\osL(f)\leq b<0$. Given any two trajectories $x(t),\bar{x}(t)$ of
a strictly contracting system, one has
$ \norm{x(t)-\bar{x}(t)}_{2,P^{1/2}} \leq e^{b(t-s)}
\norm{x(s)-\bar{x}(s)}_{2,P^{1/2}}$ for any $t\geq s\geq 0$.

The following result, derived from
Theorem~\ref{thm:explicit delta}, establishes a test for strict
contractivity of the continuous-time system in~\eqref{eq:nonl system} on
the basis of a set of measurements of $\hat{\phi}$.

\begin{theorem}\longthmtitle{Data-based test for strict
    contractivity}\label{thm:contraction}
  Given data $(Y,\Phi)$, let $\Theta=\calZ(N)$, with $N$ as in
  \eqref{eq:defN}. Let $P>0$ and $z,z^\star\in \mathbb{R}^n$. Assume
  $\Phi$ has full row rank. Then,
  \[
    (z-z^\star)^\top P \theta^\top(b(z)-b(z^\star)) \leq \gamma
    \norm{z-z^\star}^2_{2,P^{1/2}}\] for all $\theta\in\Theta$ if and
  only if
 \vspace{-0.3em} \begin{align*}
    (z-z^\star)^\top
    &P(-N_{12}N_{22}^{-1})(b(z)-b(z^\star))
    \\
    &+\sqrt{(z-z^\star)^\top
      P
      (N|N_{22})P (z-z^\star)}
    \\
    &\quad\cdot\sqrt{(b(z)-b(z^\star))^\top(-N_{22}^{-1})(b(z)-b(z^\star))}
    \\
    &\leq \gamma \norm{z-z^\star}^2_{2,P^{1/2}}.\vspace{-0.5em}
  \end{align*}
\end{theorem}
\vspace{-0.5em}As a consequence of this result, we can provide a test
to establish strict contractivity in terms only of the least-squares
estimate $\phi^{\lse}$, the data matrix $N$, and a Lipschitz constant
for the basis functions.

\begin{corollary}\longthmtitle{Strict contractivity in terms of the
    least-squares estimate}\label{cor:contraction}
  Given data $(Y,\Phi)$, let $\Theta=\calZ(N)$, with $N$ as in
  \eqref{eq:defN}. Assume $\Phi$ has full row rank.  Let $L$ be such
  that
  $ \norm{b(z)-b(z^\star)}_{2,(-N_{22})^{-1/2}} \leq L
  \norm{z-z^\star}_2$, for all $z,z^\star\in \mathbb{R}^n$.  Then, for
  $P>0$,
\vspace{-0.3em}  \[\vspace{-0.3em}
    (z-z^\star)^\top P \theta^\top(b(z)-b(z^\star)) \leq \gamma
    \norm{z-z^\star}^2_{2,P^{1/2}}
  \]
  for all $\theta\in\Theta$ and $z,z^\star\in \mathbb{R}^n$ if
  \vspace{-0.3em}\begin{equation}\label{eq:osL lse}
    \osL(\phi^{\lse}) < \gamma -
    L\sqrt{\lambda_{\textup{max}}(N|N_{22})}
    \frac{\lambda_{\textup{max}}(P)}{\lambda_{\textup{min}}(P)}.  \vspace{-0.3em}
  \end{equation}
\end{corollary}
\begin{remark}\longthmtitle{Contraction w.r.t. different norms}
  {\rm The results of this section can be readily generalized to
    contractivity (and thus one-sided Lipschitzness) with respect to
    $\norm{\cdot}_{p,R}$ for $p\in[1,\infty)$. To do this, one
    replaces the one-sided Lipschitz condition~\eqref{eq:osl
      condition} with the respective condition
    from~\cite[Table~I]{AD-SJ-FB:22} and adjust the statements
    accordingly. Note that the case of $\norm{\cdot}_{\infty,R}$ is
    not amenable to this treatment because the corresponding one-sided
    Lipschitz condition $\osL(f)$ is not linear in $f$. \oprocend}
\end{remark} 

\subsection{Suboptimal regulation of unknown systems}
\label{sec:subopt regulation}
Consider the problem of regulating an unknown linear system to a
suboptimal point of the norm of an unknown cost
function~$\hat{\phi}:\mathbb{R}^n\rightarrow\mathbb{R}^m$ on the basis
of measurements.  We consider noisy measurements~\eqref{eq:Y-W} as
described in Section~\ref{sec:problem}. In short,
$\hat{\phi}=\phi^\theta$ for some $\theta\in\Theta=\calZ(N)$.
Moreover, let
\vspace{-0.4em}\begin{equation}\label{eq:dynamics}
 \vspace{-0.4em} x_{k+1} = \hat{A}x_k+\hat{B}u_k,
\end{equation}
with state $x_k\in\mathbb{R}^n$ and input $u_k\in\mathbb{R}^r$.  Here,
$\hat{A}\in\mathbb{R}^{n\times n}$ and
$\hat{B}\in\mathbb{R}^{n\times r}$.  {\color{black} We are interested
  in regulating the system \eqref{eq:dynamics} to an equilibrium
  $x^\star$ for which we can guarantee that
  $\norm{ \hat{\phi}(x^\star)}_2\leq \delta$, with a value of $\delta$
  as small as possible.

  If $\hat{A}$ and $\hat{B}$ are known and $\hat{A}$ is stable, any
  fixed point (i.e., equilibrium) of the dynamics \eqref{eq:dynamics}
  corresponds to a constant input $u^\star$. Indeed, if we apply the
  constant input $u_k = u^\star$, the state asymptotically converges
  to the fixed point
  \vspace{-0.3em}\[ \vspace{-0.3em} \lim_{k\rightarrow \infty} x_k =
    x^\star:= (I\!-\! \hat{A})^{-1} \hat{B}u^\star,
\]
regardless of the initial condition $x_0$. This means that we can find
a static input to regulate the system towards any state
in~$\im (I\!-\! \hat{A})^{-1}\hat{B} \subseteq \mathbb{R}^n$. Hence,
if the matrices $\hat{A}$ and $\hat{B}$ are known, the problem
corresponds to Problem~\ref{prob:cautious-opt}\ref{prob:cautious} with
$\calS= \im(I\!-\! \hat{A})^{-1}\hat{B}$.
 
In the following, we assume that we do \textit{not} have access to the
matrices $\hat{A}$ and $\hat{B}$. Instead, we assume that we have not
only measurements of $\hat{\phi}$, but additionally
measurements corresponding to the unknown system. To be precise, we
collect measurements of the states and input of the system
\vspace{-0.4em}\begin{equation}\label{eq:dist lin system} 
\vspace{-0.4em}x_{k+1} = \hat{A}x_k + \hat{B} u_k +w_k.
\end{equation}
Using the basis functions
\vspace{-0.5em}\[\vspace{-0.8em}
  \psi_i(x,u) =
  \begin{cases}
    x_i & \textrm{ for } i =1,\ldots, n
    \\
    u_i & \textrm{ for } i = n+1, \ldots n+r,
  \end{cases}
\]
we can perform a linear version of set-valued regression (cf. 
Section~\ref{sec:problem}) to obtain a set $\Sigma $ such that 
$(\hat{A},\hat{B})^\top\in \Sigma = \calZ(M)$, where $M\in\mathbb{S}^{n+(n+r)}$ is
defined analogously to $N$ in~\eqref{eq:defN}. In fact, this is 
equivalent to the setup of e.g.,~\cite{HJVW-MKC-JE-HLT:22}. 

In line with all our previous reasoning, we can formalize the problem: 

\begin{problem}[Suboptimal regulation] \label{prob:subopt}
  Consider an unknown function $\hat{\phi}=\phi^\theta$ for 
  some $\theta\in\Theta=\calZ(N)$ and an unknown linear system 
  $(A,B)\in\Sigma=\calZ(M)$, where both corresponding noise models satisfy 
  Assumption~\ref{as:noise model}. Find an input $u^\star$ such that 
  \textit{each} system $(A,B)\in \Sigma$ converges to a 
  (potentially different) equilibrium which is a suboptimal point 
  for \textit{each} function $\phi^\theta$ with $\theta \in \Theta = \calZ(N)$.
\end{problem}

Our approach to this problem takes the following steps: 
\begin{itemize}
\item First, we provide conditions for the unknown system matrix
  to be stable;
\item Next, we provide bounds on the set of all equilibria resulting
  from the application of a given input to the set of
  systems~$\Sigma$;
\item Lastly, we leverage structure of these results to solve the
  suboptimal regulation problem.
\end{itemize}  }

To be able to regulate the unknown system in the same manner as
before, recall that we require that $\hat{A}$ is stable. Given that,
on the basis of measurements we can not distinguish $\hat{A}$ from
other matrices $A$ such that $(A,B)^\top\in \Sigma$, we therefore
require all such $A$ to be stable.  In the following, we describe a
test for the existence of a shared quadratic Lyapunov function, which
is a stronger stability condition. For a detailed discussion on the
conservativeness of this assumption, see
e.g.,~\cite{HJVW-MKC-HLT:22}. Now we can easily adapt
e.g.,~\cite[Thm. 5.1]{HJVW-MKC-JE-HLT:22} to obtain the following.

\begin{lemma}[Informativity for stability]\label{lem:stab}
  Let $\Sigma=\calZ(M)$, where the corresponding noise model satisfies
  Assumption~\ref{as:noise model}. Then, there exists
  $P\in\mathbb{S}^{n}$ with $P>0$ and $APA^\top<P$ for all
  $(A,B)^\top\in\Sigma$ if and only if there exist
  $\bar{P} \in\mathbb{S}^{n}$ with $\bar{P}>0$ and $\beta>0$
  such that
 \vspace{-0.4em} \begin{equation}\label{eq:LMIstab}
    {\scriptsize\begin{pmat}[{..}]
      \bar{P}-\beta I_n & 0		& 0  \cr
      0 			& -\bar{P} 	& 0  \cr
      0 			& 0 	& 0 \cr
    \end{pmat}}   -  M \geq 0 . 
  \end{equation}
\end{lemma}
%

\begin{remark}\longthmtitle{Stability and
    contractivity}\label{rem:contraction}
  {\rm Applying a fixed input $u^\star$ to the unknown system yields
    an autonomous discrete-time system. We know that for any initial
    condition this converges to a fixed point. Moreover,
    if $P>0$ and $APA^\top<P$ for all $(A,B)^\top\in\Sigma$, then
    $A^\top P^{-1} A < P^{-1}$. Therefore, the system 
    $x_{k+1} = \hat{A}x_k + \hat{B}u^\star$ 
    is strongly contracting with respect to
    $\norm{\cdot}_{2,P^{-1/2}}$.  \oprocend}
\end{remark}

When each system matrix $A$ in the set $\Sigma$ of systems consistent
with the data is stable, we can characterize the fixed points
resulting from applying the same input to each.

\begin{lemma}[Characterizing fixed points]\label{lem:fixed points}
  Let $\Sigma=\calZ(M)$, where the corresponding noise model satisfies
  Assumption~\ref{as:noise model} and $M$ has at least one positive
  eigenvalue. Assume $A$ is stable for all $(A,B)^\top \in\Sigma$.
  Given $x\in\mathbb{R}^n$ and $u\in\mathbb{R}^r$,
  $\norm{x - (I\!-\! A)^{-1}Bu}_2\leq \varepsilon$ for all
  $(A,B)^\top\in\Sigma$ if and only if there exists $\gamma\geq 0$
  such that
  \vspace{-0.4em}\begin{equation}\label{eq:LMI bd i-o}  
    \scriptsize\begin{pmat}[{..|}]
      I_n & -I_n &0 &x \cr
      -I_n & I_n &0 &-x \cr
      0& 0& 0& -u \cr\-
      x^\top & -x^\top & -u^\top & \varepsilon^2 \cr
    \end{pmat}\normalsize - 
    \gamma\begin{pmat}[{|}]
      M & 0 \cr\- 0 & 0 \cr
    \end{pmat}\geq 0	\vspace{-0.4em}
  \end{equation}
\end{lemma}
 \vspace{-1.4em}\begin{pf}
Let $(A,B)^\top \in\Sigma$. First, note that we have 
  \vspace{-0.3em}\begin{align*}
    & \qquad \norm{x - (I_n\!-\! A)^{-1}Bu}_2\leq \varepsilon
    \\  
    & \Leftrightarrow (x - (I_n\!-\! A)^{-1}Bu)(x - (I_n\!-\! A)^{-1}Bu)^\top \leq
      \varepsilon^2 I.\vspace{-0.3em}
  \end{align*}
  Since $A$ is stable, we know that $I-A$
  is nonsingular. Thus, the previous holds if and only if
  \vspace{-0.3em}\[ \vspace{-0.3em}((I\!-\! A)x-Bu)((I\!-\! A)x-Bu)^\top \! \leq \varepsilon^2
  (I\!-\! A)(I\!-\! A)^\top, \]
  or equivalently, 
  \vspace{-0.5em}\[
    \vspace{-0.3em}\scriptsize\begin{bmatrix}
      I_n
      \\
      A^\top
      \\
      B^\top
    \end{bmatrix}^{\!\top}
    \!\!
    \left(
      \begin{bmatrix}
        \varepsilon^2 I_n & -\varepsilon^2 I_n \!& 0
        \\
        -\varepsilon^2 I_n &\varepsilon^2 I_n & 0
        \\
        0 &0 &0
      \end{bmatrix}
      \! -\!
      \begin{pmatrix}
        x \\ -x \\
        -u
      \end{pmatrix} \!\!\!
      \begin{pmatrix}
        x
        \\
        -x
        \\
        -u
      \end{pmatrix}^{\!\!\top}
    \right)
    \!\!
    \begin{bmatrix}
      I_n
      \\
      A^\top
      \\
      B^\top
    \end{bmatrix}
    \geq
    0.
  \]
  The set of $(A,B)$ that satisfy the condition can be written as the
  solution set of a QMI. By assumption $\Sigma=\calZ(M)$. We can then
  apply the matrix S-Lemma~\cite[Thm 4.7]{HJVW-MKC-JE-HLT:22} and use
  a Schur complement to prove the statement.  \QED
\end{pf}
\vspace{-1.3em}Motivated by Lemma~\ref{lem:fixed points}, we define
 \vspace{-0.4em}\begin{equation}\label{eq:def eps-}
   \vspace{-0.4em}\varepsilon^-(x) := \min \{
  \varepsilon \mid \exists \gamma\geq 0, u\in\mathbb{R}^r \textrm{
    s.t. } \eqref{eq:LMI bd i-o} \textrm{ holds} \}.
\end{equation}
Given $x\in\mathbb{R}^N$, the function $\varepsilon^-(x)$ thus gives
the minimal radius of a ball around $x$, for which there exists a
$u^\star$ such that the corresponding set of fixed points of each
system in $\Sigma$ is contained in this ball. The next result
describes useful properties of this function.

\begin{lemma}\longthmtitle{Properties of the minimal radius function}
  Let $\Sigma=\calZ(M)$, where the corresponding noise model satisfies
  Assumption~\ref{as:noise model} and $M$ has at least one positive
  eigenvalue. Assume $A$ is stable for all $(A,B)^\top \in\Sigma$.
  Then $\varepsilon^-(x) \leq \norm{x}_2$ for any $x\in\mathbb{R}^n$
  and $\varepsilon^-$ is convex.
\end{lemma}
 \vspace{-1em}\begin{pf}
  The statement $\varepsilon^-(x) \leq \norm{x}_2$ follows by noting
  that, for any $x$ and $M$, \eqref{eq:LMI bd i-o} holds with $u=0$,
  $\gamma=0$ and $\varepsilon=\norm{x}_2$. Next, suppose that for all
  $(A,B)^\top\in\Sigma$,
 \vspace{-0.4em}  \[
    \norm{x - (I\!-\! A)^{-1}Bu}_2\leq \varepsilon, \quad \norm{\bar{x}
      - (I\!-\! A)^{-1}B\bar{u}}_2\leq \bar{\varepsilon}.
  \]
  Let $\lambda\in [0,1]$ and define
  $x_\lambda:= \lambda x + (1-\lambda)\bar{x}$ and similarly for
  $u_\lambda$ and $\varepsilon_\lambda$. Then,
  \vspace{-0.4em} \begin{align*}
    & \norm{x_\lambda - (I\!-\! A)^{-1}Bu_\lambda}_2
      \leq \lambda\norm{x - (I\!-\! A)^{-1}Bu}_2
    \\
    & +(1-\lambda)\norm{\bar{x}
      - (I\!-\! A)^{-1}B\bar{u}}_2
      \leq \varepsilon_\lambda.
  \end{align*}
  Therefore, $\varepsilon^-(x_\lambda)\leq \varepsilon_\lambda$,
  proving that $\varepsilon^-$ is convex in~$x$. \QED
\end{pf}
 \vspace{-0.6em}{\color{black}The pieces are now in place for our solution to 
	Problem~\ref{prob:subopt} in the following result. }

\begin{theorem}\longthmtitle{Suboptimal
    regulation}\label{thm:subopt reg}
  Let $\Sigma=\calZ(M)$ and $\Theta=\calZ(N)$, where the corresponding
  noise models satisfy Assumption~\ref{as:noise model}. Assume $A$ is
  stable for all $(A,B)^\top \in\Sigma$ and $\phi^\theta$ is Lipschitz
  with constant $L$ for all $\theta\in\Theta$. Then
  \begin{equation}\label{eq:subopt reg}
    \min_u \norm{\hat{\phi}((I-\hat{A})^{-1}\hat{B}u)}_2 \leq \min_x \,
    \big(g(x;\Theta)+L\varepsilon^-(x)\big).
  \end{equation}
\end{theorem}
\begin{pf}
  Using the triangle inequality and Lipschitzness,
  \vspace{-0.4em} \[ 
    \norm{\hat{\phi}((I-\hat{A})^{-1}\hat{B}u)}_2 \leq
    \norm{\hat{\phi}(x)}_2 + L\norm{x-(I\!-\!
      \hat{A})^{-1}\hat{B}u}_2,
  \]
  for any $x \in \mathbb{R}^n$. Note that this separates the unknown
  parameter $\hat{\theta}$ from the unknown pair
  $(\hat{A},\hat{B})$. From here, we deduce
  \vspace{-0.5em}\begin{align*}
    & \min_u \norm{\hat{\phi}((I-\hat{A})^{-1}\hat{B}u)}_2
    \\
    & \leq \min_{x,u} \Big(\max_{\theta\in\Theta}
      \norm{\phi^\theta(x)} +
      L\Big(\max_{(A,B)\in\Sigma}\norm{x-(I\!-\!
      A)^{-1}Bu} \Big)\Big).
  \end{align*}
  Since only the second term is dependent on $u$, we can use the
  definition of $\varepsilon^-(x)$ and Lemma~\ref{lem:fixed points} to
  conclude the statement.  \QED
\end{pf}
 \vspace{-1.2em}{\color{black}Similar to Problem~\ref{prob:cautious-opt}\ref{prob:cautious}, we 
can now resolve the optimization problem on the right hand side under 
various regularity conditions. To illustrate this, recall that $\varepsilon^-$ 
is convex, and thus we can resolve this
problem efficiently if so is $g(\cdot;\Theta)$ (see
Corollary~\ref{cor:conv uncert}).} We conclude this section by noting
that Theorem~\ref{thm:subopt reg} requires only stability of each
system matrix $A$. Under the stronger assumption of quadratic
stability, which can be determined using Lemma~\ref{lem:stab}, we can
provide a transient guarantee on the values of the unknown function.

\begin{theorem}\longthmtitle{Transient values of the cost function}
	\label{thm:transients}
  Let $\Sigma=\calZ(M)$ and $\Theta=\calZ(N)$, where the corresponding
  noise models satisfy Assumption~\ref{as:noise model}. Let $P>0$ such
  that $APA^\top<P$ for all $(A,B)^\top \in\Sigma$ and assume
  $\phi^\theta$ is Lipschitz with constant $L$ for all
  $\theta\in\Theta$. Let $\delta>0$ and $u^\star$ be a fixed input
  such that
  $\norm{\hat{\phi}((I-\hat{A})^{-1}\hat{B}u^\star)}_2 \leq \delta$.
  Given an initial condition $\hat{x}_0\in\mathbb{R}^n$, let
  $\hat{x}_k$ denote the trajectory of the unknown system
  $(\hat{A},\hat{B})$ corresponding to the input $u^\star$ and
  starting from $\hat{x}_0$. Then, there exists $\lambda\in[0,1)$ such
  that
  \vspace{-0.5em}\[
   \vspace{-0.3em} \norm{\hat{\phi}(\hat{x}_k)}_2 \leq \delta +
    \lambda^kL\frac{\lambda_{\textup{max}}(P)}{\lambda_{\textup{min}}(P)}
    \norm{\hat{x}_{0}-x^\star}_2.
  \]
\end{theorem}
 \vspace{-1.4em}\begin{pf}
  Let $x^\star$ be the fixed point corresponding to the
  input~$u^\star$ for the system $(\hat{A},\hat{B})$.
  If $APA^\top<P$ for all $(A,B)^\top \in\Sigma$, then the set 
  $\{A\mid (A,B)^\top\in\Sigma\}$ is bounded. Since the set is closed
  by definition, we can conclude that $APA^\top<P$ for all 
  $(A,B)^\top \in\Sigma$ if and only if there exists $\lambda\in[0,1)$ 
  such that $APA^\top <\lambda P$ for all $(A,B)^\top \in\Sigma$.
  In particular, this implies that the \textit{affine} system
  resulting from the application of input $u^\star$ to
  $(\hat{A},\hat{B})$ is strongly contracting
  (cf. Remark~\ref{rem:contraction}). Thus, for any $k\geq1$,
 \vspace{-0.4em}  \[ 
    \norm{\hat{x}_{k}-x^\star}_{2,P^{-1/2}} \leq \lambda
    \norm{\hat{x}_{k-1}-x^\star}_{2,P^{-1/2}}.
  \]
  By applying this repeatedly, we conclude
  $ \norm{\hat{x}_{k}-x^\star}_{2,P^{-1/2}} \leq
  \lambda^{k}\norm{\hat{x}_0-x^\star}_{2,P^{-1/2}}$, and thus
  \vspace{-0.5em}\[
    \vspace{-0.3em}\norm{\hat{x}_{k}-x^\star}_2 \leq
    \lambda^{k}\frac{\lambda_{\textup{max}}(P^{-1})}{\lambda_{\textup{min}}(P^{-1})}
    \norm{\hat{x}_0-x^\star}_2.
  \]
  Using the triangle inequality and the fact that $\hat{\phi}$ is
  Lipschitz with constant $L$,
 \vspace{-0.4em}  \[
    \norm{\hat{\phi}(\hat{x}_k)}_2\leq
    \norm{\hat{\phi}(x^\star)}_2+L\norm{\hat{x}_{k}-x^\star}_2.
  \] 
  Combining this previous with the fact that
  $ \lambda_{\textup{max}}(P^{-1}) = {1}/{\lambda_{\textup{min}}(P)}$
  and
  $ \lambda_{\textup{min}}(P^{-1}) = {1}/{\lambda_{\textup{max}}(P)}$
  yields
  \vspace{-0.3em}\[
  \vspace{-0.3em}  \norm{\hat{\phi}(\hat{x}_k)}_2 \leq \delta + \lambda^{k}
    L\frac{\lambda_{\textup{max}}(P)}{\lambda_{\textup{min}}(P)}
    \norm{\hat{x}_{0}-x^\star}_2 .
  \]
  %
  This proves the statement.\QED
\end{pf}
\begin{remark}\longthmtitle{Suboptimal regulation of unstable systems}
  {\rm Our discussion above assumes that the unknown system and all of
    those consistent with the measurements are stable. When this is
    not the case, we can instead proceed by
    first testing if the data is sufficiently informative for
    \emph{quadratic stabilization}, i.e., to guarantee the existence
    of a static feedback $K$ such that $(A+BK)P(A+BK)^\top<P$ for all
    $(A,B)^\top\in \Sigma$, cf~\cite[Thm. 5.1]{HJVW-MKC-JE-HLT:22}.
    Then, we can characterize the behavior of the fixed points of the
    resulting systems by considering
    \vspace{-0.3em}\[
     \vspace{-0.3em} \norm{x - (I\!-\! A+BK)^{-1}Bu}_2\leq \varepsilon.
    \]
    Adapting Lemma~\ref{lem:fixed points} and the definition of
    $\varepsilon^-$ accordingly is then straightforward. This gives
    rise to a two-step approach: (i) determine a gain $K$ which
    stabilizes all systems in $\Sigma$ and (ii) minimize
    $g(x;\Theta)+L\varepsilon^-(x)$. Note that the choice of $K$ among
    potentially many stabilizing gains influences the resulting
    $\varepsilon^-$ in a nontrivial way. We leave the analysis of
    jointly minimizing $g(x;\Theta)+L\varepsilon^-(x)$ over all $x$
    and all possible choices $K$ as an open problem. \oprocend}
\end{remark}


 \vspace{-0.6em}\section{Online cautious optimization}\label{sec:online}

Our exposition so far has considered one-shot optimization settings
where, given a set of measurements, we determine regularity
properties and address suboptimization of the unknown function. In
this section, we consider \textit{online} scenarios, where one can repeatedly
combine optimization of the unknown function with the collection
of additional measurements, aiming at reducing the optimality
gap. Consistent with our exposition, this means that we are interested
in finding guaranteed upper bounds on $c^\top \hat{\phi}$ on the basis
of initial measurements, then refining this bound on the basis of
newly collected measurements\footnote{To keep the
  presentation contained, we do not consider minimization of
  $\norm{\hat{\phi}}_2$, but the results can be adapted to this
  case.}. We aim at doing this in a monotonic fashion, that
is, in a way that has the obtained upper bounds do not increase as
time progresses.

Consider repeated measurements of the form
 \vspace{-0.4em}\begin{equation}\label{eq:measure} \vspace{-0.4em}
  Y_k = \hat{\theta}^\top\Phi_k+ W_k, \textrm{ with } W_k^\top \in
  \calZ(\Pi), 
\end{equation}
where $Y_k$ and $W_k$ are as in~\eqref{eq:Y-W} and $k\in\mathbb{N}$. Let
$ \Theta_k:= \calZ(N_k)$ be the set of of parameters which are
compatible with the $k^{\textup{th}}$ set of measurements, where
\vspace{-0.5em}\[
\vspace{-0.3em}  N_k:= {\scriptsize\begin{bmatrix} I & Y_k \\ 0& -\Phi_k \end{bmatrix}}
  \Pi {\scriptsize\begin{bmatrix} I & Y_k \\ 0& -\Phi_k \end{bmatrix}}^\top.
\]
Our online optimization strategy is described as follows.

\begin{algo}\longthmtitle{Online cautious
    optimization}\label{alg:onl caut subopt}
  Consider an initial candidate optimizer $z_0\in\mathbb{R}^n$, an
  initial set of data $(Y_0,\Phi_0)$, and $c\in\mathbb{R}^m$. Define
  the initial set of compatible parameters as
  $\bar{\Theta}_0 = \calZ(N_0)$, where $N_0$ is given as
  in~\eqref{eq:defN}. For $k\geq 1$, we alternate two steps:
  \begin{enumerate}
  \item Employ the set $\bar{\Theta}_{k-1}$ to improve the candidate
    optimizer by finding $z_k$ such that
 \vspace{-0.3em}   \begin{equation}\label{eq:z decrease}
      g_c(z_k;\bar{\Theta}_{k-1}) \leq
      g_c(z_{k-1};\bar{\Theta}_{k-1}).
    \end{equation}
  \item Measure the function as in \eqref{eq:measure} and update the
    set of parameters by finding a set $\bar{\Theta}_k$
    such that $\hat{\theta}\in\bar{\Theta}_k$ and
\vspace{-0.3em}    \begin{equation}\label{eq:Theta decrease}
      g_c(z_k;\bar{\Theta}_{k}) \leq g_c(z_{k};\bar{\Theta}_{k-1}).
    \end{equation}
  \end{enumerate} 	
\end{algo}

In the following, we enforce~\eqref{eq:Theta decrease} by considering
parameters consistent with \textit{all} previous measurements,
\vspace{-0.3em}\begin{equation}\label{eq:update theta}
  \vspace{-0.3em}\bar{\Theta}_k := \Theta_0\cap\ldots\cap\Theta_{k}.
\end{equation}
Note that Algorithm~\ref{alg:onl caut subopt} provides a sequence
of upper bounds to true function values,
on the basis of measurements.  This means that after \textit{any}
number of iterations, we obtain a worst-case estimate of the function
value $c^\top\hat{\phi}(z_k)$ and, as such, for the minimum of
$c^\top\hat{\phi}$. To be precise, for each $k\geq 1$,
\vspace{-0.3em}\begin{equation}\label{eq:ubounds alg}
\vspace{-0.3em}  \min_{z\in\mathbb{R}^n} c^\top\hat{\phi}(z)\leq
  c^\top\hat{\phi}(z_k)\leq  g_c(z_k;\bar{\Theta}_{k-1}), 
\end{equation}
and these bounds are monotonically nonincreasing
\vspace{-0.3em}\begin{equation}\label{eq:monotone}
\vspace{-0.3em}  g_c(z_{k+1};\bar{\Theta}_{k}) \leq g_c(z_k;\bar{\Theta}_{k-1})  
\end{equation}
Therefore this sequence of upper bounds is nonincreasing and bounded
below, and we can conclude that the algorithm converges. However,
without further assumptions, one cannot guarantee convergence of the
upper bounds to the minimal value of $c^\top\hat{\phi}$, or
respectively of $z_k$ to the global minimum
of~$c^\top\hat{\phi}$. This is because, even though when repeatedly
collecting measurements, it is reasonable to assume that the set of
consistent parameters would decrease in size, this is not necessarily
the case in general. In particular, a situation might arise where
repeated measurements corresponding to a worst-case, or adversarial,
noise signal give rise to convergence to a fixed bound with nonzero
uncertainty. {\color{black} This is the point we address next by
  considering random, stochastic noise realizations}.

\subsection{Collection of uniformly distributed measurements}
{\color{black} Recall that, so far we have assumed that the noise, as
  collected in the matrix $W$, satisfies Assumption~\ref{as:noise
    model}, i.e., we have access to $\Pi\in\mathbb{S}^{m+T}$ with
  $\Pi_{22}<0$ and $\Pi|\Pi_{22} \geq 0$, such that
  $W^\top \in \calZ(\Pi)$. In particular, this implies that
  $\calZ(\Pi)$ is bounded. In this section, we consider a scenario
  where the set $\calZ(\Pi)$ has nonempty interior (equivalent to
  $\Pi|\Pi_{22} >0$,~\cite[Thm. 3.2]{HJVW-MKC-JE-HLT:22}) and the
  noise samples are not only bounded but \textit{distributed uniformly
    randomly} over the set $\calZ(\Pi)$. To formalize this, consider a
  measure $\mu$ on $\mathbb{R}^{m\times T}$, and define a probability
  density function $p$ by
\[
  p(W) :=
  \begin{cases}
    \frac{1}{\mu(\calZ(\Pi))} & W^\top\in\calZ(\Pi),
    \\ 0 &
           \textrm{otherwise}.
  \end{cases}
\]
As a notational shorthand, we write $W^\top\sim\uniform(\calZ(\Pi))$
if the distribution of $W$ follows the probability density function
$p$.}  The following result shows that, under such uniformly
distributed noise samples, the size of the set of parameters indeed
shrinks with increasing $k$.

\begin{theorem}\longthmtitle{Repeated measurements leading to
    shrinking set of consistent parameters}\label{thm:repeats}
  Suppose that the measurements $(Y_k,\Phi_k)$ are collected such that
  $W_k^\top\sim \uniform(\calZ(\Pi))$ and $\sigma_-(\Phi_k(z))\geq a$
  for some $a\in\mathbb{R}_{>0}$ and all $z\in \mathbb{R}^n$.  Then,
  for any probability $0<\pi<1$ and any $\varepsilon>0$, there exists
  $k$ such that if $B_\varepsilon(\hat{\theta})$ denotes the ball of
  radius $\varepsilon>0$ centered at $\hat{\theta}$ in
  $\mathbb{R}^{k\times m}$, then \vspace{-0.3em}\begin{equation}
    \label{eq:delta-probability}
    \vspace{-0.3em}p(\Theta_0\cap\ldots \cap \Theta_{k-1}\subseteq
    B_\varepsilon(\hat{\theta}))>\pi.   
  \end{equation}
\end{theorem}
 \vspace{-1.4em}\begin{pf}
  Let $k\geq 0$ and consider the set $\Theta_k=\calZ(N_k)$.  In terms
  of the true parameter $\hat{\theta}$, we obtain
  $Y_k = \hat{\theta}^\top\Phi_k +W_k$, with
  $W_k^\top\sim \uniform(\calZ(\Pi))$. This implies that
  \vspace{-0.3em}\[ \vspace{-0.3em} Y_k^\top \sim \uniform (\calZ(\Pi)
    +\Phi_k^\top \hat{\theta}).
  \]
  For $\theta\in\mathbb{R}^k$, we know that $\theta\in\Theta_k$ if and
  only if $Y_k^\top -\Phi_k^\top\theta\in \calZ(\Pi)$. Therefore,
  $ Y_k^\top-\Phi_k^\top\theta \sim \uniform (\calZ(\Pi) +\Phi_k^\top
  (\hat{\theta}-\theta))$.  We can now calculate the probability of
  $\theta\in\Theta_k$ to be
  \vspace{-0.3em}\[
    p(\theta\in\Theta_k) =
    \frac{\mu((\Phi_k^\top(\hat{\theta}-\theta)+\calZ(\Pi))\cap
      \calZ(\Pi))}{\mu(\calZ(\Pi))}.\vspace{-0.3em}
  \]
  Recall, however, that the true parameter $\hat{\theta}$ is
  unknown. Let $\theta\not\in B_\varepsilon(\hat{\theta})$, then,
  since by assumption $\sigma_-(\Phi_k)\geq a$, we can derive that
 \vspace{-0.3em}\[
    p(\theta\in\Theta_k)
     \leq \max_{\norm{V}_2=1} \frac{\mu((a\varepsilon
        V+\calZ(\Pi))\cap \calZ(\Pi))}{\mu(\calZ(\Pi))}=:
        \eta(a\varepsilon)  .\vspace{-0.3em}
  \]
  Note that, since $\calZ(\Pi)$ is bounded and convex, we have in
  addition that $\eta(\lambda)<1$ for $\lambda>0$. This means that,
  regardless of the values of the measurements, the probability that
  $\theta\in\Theta_k$ is strictly less than 1. Now, by repeating such
  measurements, we see that for
  $\theta\not\in B_\varepsilon(\hat{\theta})$, we have
  \vspace{-0.3em}\[
    p(\theta\in \Theta_0\cap \ldots \cap \Theta_{k-1}) =
    p(\theta\in\Theta_0)\cdots p(\theta\in\Theta_{k-1}) \leq
    \eta(a\varepsilon)^k .\vspace{-0.3em}
  \]
  From this, if $k$ is large enough, then~\eqref{eq:delta-probability} 
  holds.  \QED
\end{pf}

 \vspace{-1.4em}\begin{remark}\longthmtitle{Extension to other probability
    distributions} {\rm Note that the proof of
    Theorem~\ref{thm:repeats} essentially depends on the fact that
    $p(W)\neq 0 $ for all $W^\top \in \calZ(\Pi)$. Therefore, we can
    draw the same conclusion for other probability distributions that
    satisfy this assumption, such as truncated Gaussian
    distributions. \oprocend }
\end{remark}

Thus, if the noise is uniformly random, Theorem~\ref{thm:repeats}
guarantees that the set of parameters converges. This allows us to
conclude that, under Algorithm~\ref{alg:onl caut subopt}, the sequence
$ g_c(z_k;\bar{\Theta}_k)$ converges to the sequence
$c^\top\hat{\phi}(z_k)$ pointwise.

\subsection{Methods to improve the candidate optimizer}

In this section we discuss methods to guarantee that \eqref{eq:z
  decrease} holds.  A direct way of enforcing this is by simply
optimizing the expression, that is, updating $z_k$ as
{\color{black} \vspace{-0.4em}\begin{equation}\label{eq:global opt} \vspace{-0.4em}
  z_k \in \argmin_{z\in \mathbb{R}^n}
  g_c(z;\bar{\Theta}_{k-1}),\vspace{-0.3em}
\end{equation}
where  $\argmin$ denotes the set of arguments of the minimization.}
This requires us to find a global optimizer of the function
$g_c(z;\bar{\Theta}_{k-1})$, which might be difficult to obtain
depending on the scenario.
As an alternative, we devise a procedure where the
property~\eqref{eq:z decrease} is guaranteed by \textit{locally}
updating the candidate optimizer along with collecting new
measurements \textit{near} the candidate optimizer. We show that,
using only local information in this manner, the algorithm converges
to the true optimizer under appropriate technical assumptions.

In order to formalize this notion of `local', we both measure and
optimize on a polyhedral set around the candidate optimizer. For this,
let $\calF =\{f_i\}_{i=1}^T \subseteq\mathbb{R}^n$ be a finite set
such that $0\in\interior (\conv \calF)$. Let
$\calS(z) := z +\conv\calF$. Then,~\eqref{eq:z decrease} holds if we update
the optimizer by
 \vspace{-0.4em}\begin{equation}\label{eq:update est}
  \vspace{-0.4em} z_k \in  \argmin_{z\in
    \calS(z_{k-1})} g_c(z;\Theta_0\cap\ldots\cap\Theta_{k-1}).
\end{equation} 
We measure the function at all points in $z_k +\calF$, that is,
we take $\Phi_k = \Phi^\calF(z_k)$, where
 \vspace{-0.4em}\[ \vspace{-0.4em}
\Phi^\calF(z) :={\scriptsize\begin{bmatrix}
  \phi_1(z+f_1) & \dots & \phi_1(z+f_T)
  \\
  \vdots & &
             \vdots
  \\
  \phi_k(z+f_1) & \dots & \phi_k(z+f_T)
\end{bmatrix}}.
\]
The online optimization procedure then incrementally incorporates
these measurements to refine the computation of the candidate
optimizer. The following result investigates the properties of the
resulting \textit{online local descent}.

\begin{theorem}[Online local descent]\label{thm:online}
  Let $\calF$ be a finite set such that $0\in\interior (\conv
  \calF)$. Let $\calS(z) := z +\conv\calF$. Consider
  Algorithm~\ref{alg:onl caut subopt}, employing local update
  rules~\eqref{eq:update est} and~\eqref{eq:update theta}, where
  $\theta_k=\calZ(N_k)$ and $\Phi_k= \Phi^\calF(z_k)$.  Suppose that
  the initial point $z_0\in\mathbb{R}^n$ and the data
  $(Y_0,\Phi^\calF(z_0))$ are such that $c^\top\phi^\theta$ is
  strictly convex for all $\theta\in\Theta_0$.  Then the following
  hold:
  \begin{enumerate}
  \item\label{item:convex} For any $k\geq 1$, the problem
    \eqref{eq:update est} is strictly convex;
  \item\label{item:strict mono} If $z_k\neq z_{k+1}$, then
  \vspace{-0.4em}   \[ \vspace{-0.4em}g_c(z_{k+1};\bar{\Theta}_{k}) < g_c(z_k;\bar{\Theta}_{k-1}), \]
    that is, \eqref{eq:monotone} holds with a strict inequality.
  \end{enumerate}
\end{theorem}
 \vspace{-1.4em}\begin{pf}
  Statement~\ref{item:convex} follows from Proposition~\ref{cor:convexity}. 
  Then, note that if $z_k\neq z_{k+1}$, we have that
 \vspace{-0.4em}  \[ \vspace{-0.4em}
    g_c(z_{k+1};\Theta_0\cap\ldots\cap\Theta_{k})<
    g_c(z_k;\Theta_0\cap\ldots\cap\Theta_{k}),
  \]
  thus proving~\ref{item:strict mono}.\QED
\end{pf}

 \vspace{-1.3em}\begin{remark}\longthmtitle{Reduction of the complexity} {\rm Most
    optimization schemes require a large number of function
    evaluations in order to find an optimizer of
    e.g.,~\eqref{eq:update est}.  In turn, finding values of
    $g_c(\cdot;\bar{\Theta}_k)$ requires us to resolve another
    optimization problem. This nested nature of the problem means that
    reducing the number of evaluations can speed up computation
    significantly. This can be done by relaxing~\eqref{eq:update est}
    as follows.  Assume a Lipschitz constant of $c^\top\hat{\phi}$ is
    known or obtained from measurements. Since $\calF$ is a finite
    set, there exists $\nu\in\mathbb{R}_{\geq 0}$ such that each
    $z\in \calS(z_{k-1})$ is at most a distance $\nu$ from a point in
    $z_{k-1}+\calF$. As such, we can instead find
    \vspace{-0.4em} \[ \vspace{-0.4em}
      z_k \in \argmin_{z\in z_{k-1}+\calF}
      g_c(z;\Theta_0\cap\ldots\cap\Theta_{k-1}),
    \]
    and employ the value of $\nu$ and the Lipschitz constant to
    approximate the optimal value of \eqref{eq:update est}.  \oprocend
  }
\end{remark}

Recall that, without making further assumptions, we can only guarantee
that the algorithm converges, but not that it converges to the true
optimizer. Moreover, even when the assumptions of
Theorem~\ref{thm:repeats} hold, we do not yet have a way of
concluding that the algorithm has converged. 
The following result employs the uncertainty near the optimizer,
to provide a criterion for convergence.

\begin{lemma}[Stopping criterion]\label{lem:stop}
  Let $\Theta$ be compact and $\calS\subseteq\mathbb{R}^n$
  closed. Define
  \vspace{-0.3em}\[
  \vspace{-0.3em}  \bar{z} \in \argmin_{z\in\calS} g_c(z;\Theta), \quad
    \hat{z} \in \argmin_{z\in\calS} c^\top\hat{\phi}(z).
  \]
  Then
  $ g_c(\bar{z};\Theta) \geq c^\top\hat{\phi}(\hat{z})\geq
  g_c(\bar{z};\Theta) - 2\max_{z\in\calS} U_c(z;\Theta)$.
\end{lemma}
 \vspace{-0.8em}\begin{pf}
  For the first inequality, note that by definition,
  \vspace{-0.3em}\[
  \vspace{-0.3em}  g_c(\bar{z};\Theta) \geq c^\top\hat{\phi}(\bar{z}) \geq
    c^\top\hat{\phi}(\hat{z}).
  \]
  To show the second inequality, note that
  $g_c(\bar{z};\Theta)\leq g_c(\hat{z};\Theta)$, and
  \vspace{-0.3em}\[
   \vspace{-0.3em} 2\max_{z\in\calS} U_c(z;\Theta) \geq g_c(\hat{z};\Theta)
    -c^\top\hat{\phi}(\hat{z};\Theta).
  \]
  Combining these, we obtain
  \vspace{-0.3em}\[
  \vspace{-0.3em}  g_c(\bar{z};\Theta)- 2\max_{z\in\calS} U(z;\Theta)\leq
    c^\top\hat{\phi}(\hat{z};\Theta),
  \]
  which proves the result.\QED
\end{pf}
 \vspace{-1.2em}Importantly, we note that the upper and lower bounds in
Lemma~\ref{lem:stop} can be computed purely in terms of data.

If the maximal uncertainty on the set $\calS$ is equal to zero, then
the local minima of $g_c$ and $c^\top\hat{\phi}$ coincide. In
addition, if $c^\top\hat{\phi}$ is strictly convex, we have that any
local minimum in the \emph{interior} of $\calS$ is a global minimum.

\begin{corollary}\longthmtitle{Uncertainty under
    repetition}\label{coro:stochastics}
  Under the assumptions of Theorem~\ref{thm:online}, suppose in
  addition that for $k\geq 1$, the measurements
  $(Y_k,\Phi^\calF(z_k))$ are collected such that
  $W_k^\top\sim \uniform(\calZ(\Pi))$ and
  $\sigma_-(\Phi^\calF(z_k))\geq a$ for some $a\in\mathbb{R}_{>0}$ and
  all $k$. Then, for any $z\in\mathbb{R}^n$, the expected value of the
  uncertainty monotonically converges to 0, that is,
 \vspace{-0.4em}  \[ \vspace{-0.4em}
    U_c(z;\Theta_0\cap\ldots\cap \Theta_{k-1})\geq
    U_c(z;\Theta_0\cap\ldots\cap \Theta_k),
  \]
  and
  $\lim_{k\rightarrow \infty}
  \mathbb{E}(U_c(z;\Theta_0\cap\ldots\cap\Theta_k)) = 0$.
\end{corollary}
 \vspace{-0.8em}\begin{pf}
  The monotonic decrease of the uncertainty readily follows from its
  definition.  If
  $\Theta_0\cap\ldots \cap \Theta_{k-1}\subseteq
  B_\varepsilon(\hat{\theta})$, then it is immediate from the
  definitions that
 \vspace{-0.4em}  \[ \vspace{-0.4em}
    U_c(z;\Theta_0\cap\ldots \cap \Theta_{k-1}) \leq
    U_c(z;B_\varepsilon(\hat{\theta})).
  \]
  Given that the map $\theta \mapsto c^\top\phi^\theta(z)$ is linear
  in $\theta$, the right-hand side can be seen to converge to 0 for
  $\epsilon\rightarrow0$. Combining these pieces with
  Theorem~\ref{thm:repeats} proves the statement. \QED
\end{pf}
 \vspace{-0.9em}As a consequence of Lemma~\ref{lem:stop} and
Corollary~\ref{coro:stochastics}, we conclude that the expected
difference between the optimal value
$\min_{z\in\calS} c^\top\hat{\phi}(z)$ of the unknown function and the
optimal value
$\min_{z\in\calS} g_c(z;\Theta_0 \cap \dots \cap \Theta_k)$ provided
by online local descent both converge to zero.

 \vspace{-0.6em}\section{Simulation examples}\label{sec:sims}
Here we provide two simulation examples to illustrate the proposed
framework and our results. We consider an application to data-based
contraction analysis based on Section~\ref{sec:contraction}, and the
data-based suboptimal regulation of an unmanned aerial vehicle, using
Section~\ref{sec:subopt regulation}. We refer the interested reader to
our conference paper~\cite[Section~\ref{sec:online}]{JE-JC:23-cdc} for
an illustration of a simplified version of online cautious
suboptimization.

\subsection{Application to data-based contraction analysis}
Consider continuous-time systems $\dot{z}=\phi^\theta(z)$, where
$z\in\mathbb{R}^2$, the parameter $\theta\in\mathbb{R}^{6\times 2}$,
and with basis functions
\vspace{-0.3em}\begin{align*}
  \phi_1(z) = z_1, \hspace{1em}
  & \phi_3(z) = \sin(z_1)-1,
  & \phi_5(z)
    =
    \cos(z_1),
  \\ 
  \phi_2(z) = z_2, \hspace{1em}
  & \phi_4(z) = \sin(z_2)-1,
  & \phi_6(z)
    =
    \cos(z_2). \vspace{-0.3em}
\end{align*}
These define the vector-function
$b:\mathbb{R}^2\rightarrow\mathbb{R}^6$, cf.~\eqref{eq:b}. Suppose
that the parameter of the unknown function $\hat{\phi}$ is
\vspace{-0.3em}\[
  \vspace{-0.3em}\hat{\theta}^\top =
  \scriptsize\begin{bmatrix}
    -6
    & \hphantom{-}1
    &
      -1
    &
      \hphantom{-}1 & -1 & \hphantom{-}1
    \\
    \hphantom{-}0 & -6
    &
      \hphantom{-}1& -1 &\hphantom{-}1&-1
  \end{bmatrix}.
\]
We are interested in determining whether the system
$\dot{z} =\hat{\phi}(z)$ is strictly contracting with respect to the
unweighted $2-$norm on the basis of noisy measurements.

As a first step, note that {\color{black} for the Jacobian of $b$, we
  have }$\norm{J(b)}_2 = \sqrt{2}$. Then writing
$\hat{\phi}(z) = -6z +h(z)$, we can conclude that
\vspace{-0.3em}\begin{equation}\label{eq:osL phi hat}
  \begin{aligned}
    \osL(\hat{\phi})
    &\leq -6+ \osL(h)
      \leq -6 +\norm{J(h)}_2
    \\
    & \leq -6+ \norm{\left\lbrack
      \begin{smallmatrix}
        0 & \hphantom{-}1
        & -1& \hphantom{-}1 & -1 & \hphantom{-}1
        \\
        0 &
            \hphantom{-}0 & \hphantom{-}1& -1
                            &\hphantom{-}1&-1
      \end{smallmatrix}
      \right\rbrack}_2\cdot
      \norm{J(b)}_2
    \\
    & \approx -1.8694.
  \end{aligned}\vspace{-0.3em}
\end{equation}
Thus, the system $\dot{z}=\hat{\phi}(z)$ is indeed strictly
contracting.

\begin{figure}[ht]
  \centering
  \includegraphics[width=\linewidth,trim=3.5cm 2.8cm 3.5cm 1cm,
  clip]{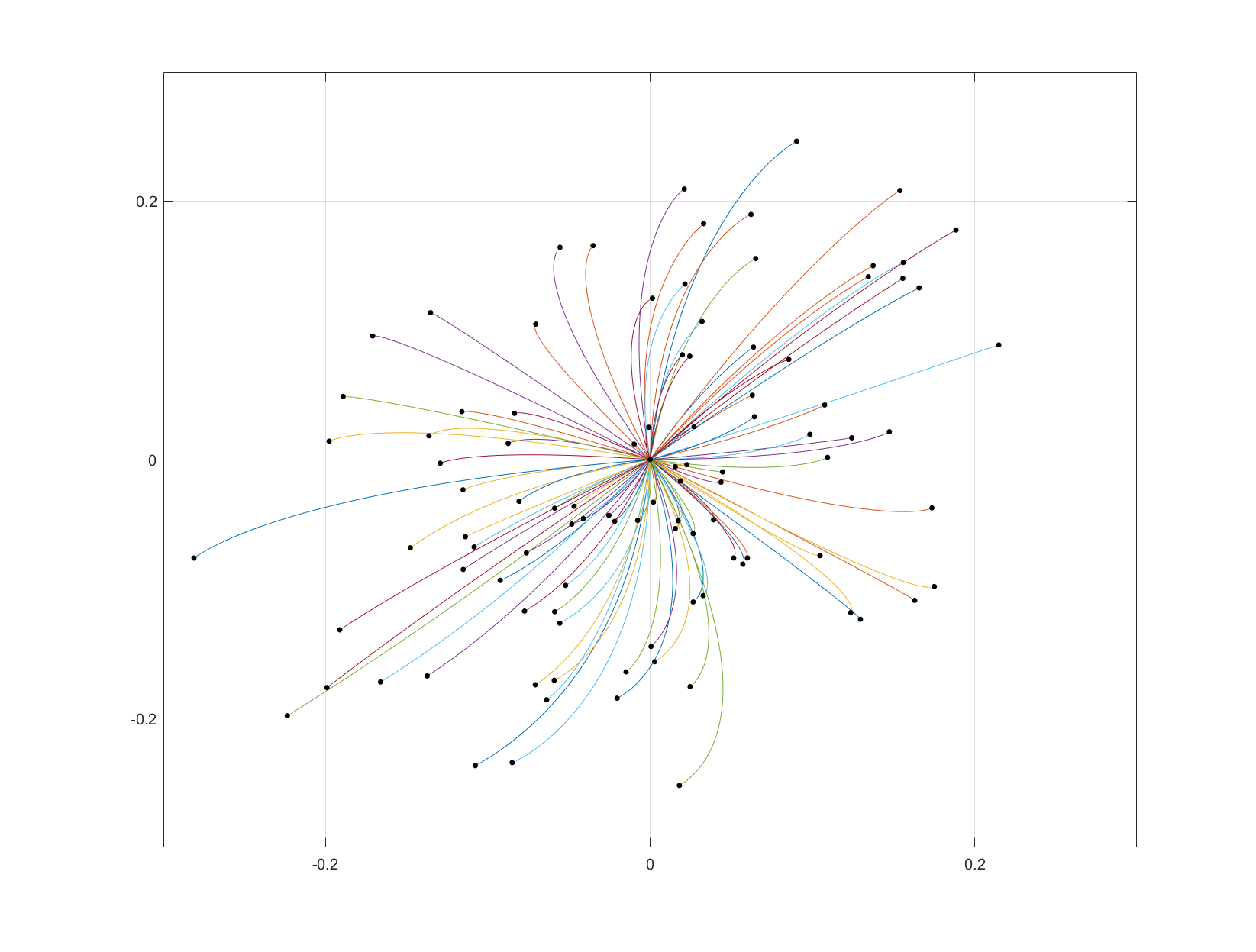}%
  \caption{Trajectories of $\dot{z}=\phi^\theta(z)$ from the origin
    for 100 randomly generated realizations of the parameter
    $\theta\in\Theta$. The black dots indicate the resulting equilibria. }\label{fig:plots systems from origin}
\end{figure}

However, as before we assume that we do not have access to the value
of $\hat{\theta}$ and can only assess contractivity based on
measurements of $\hat{\phi}$. To be precise, we assume that we have
access to samples of a single continuous-time trajectory starting from
$z(0)=(10,-20)^\top$ at the time instances
\[
  t \in \{0,0.01,\ldots,0.25\}.
\]
The noise is assumed to satisfy $WW^\top \leq 10\cdot I
_2$, giving rise to a set $\Theta=\calZ(N)$. Two examples of
parameters compatible with the measurements are given by
\begin{align*} 
  &\scriptsize\begin{bmatrix}
     -5.2588 & \hphantom{-}1.3807 & -1.4810&
                                             \hphantom{-}0.8087 & -0.8825 & \hphantom{-}1.5444
     \\  
     -0.2166 & -6.1490 & \hphantom{-}1.1110 & -0.5087
                                                                &\hphantom{-}0.4157&-0.6826\end{bmatrix},
  \\
  &\scriptsize\begin{bmatrix} -6.3111 & \hphantom{-}0.8091 & -0.6433&
                                                           \hphantom{-}1.0725 & -0.5355 & \hphantom{-}0.9938
     \\ 
     \hphantom{-}0.7750 & -5.6377 & \hphantom{-}0.9785 & -1.1375 &\hphantom{-}1.3976 &-0.8905
   \end{bmatrix}.
\end{align*}
To illustrate the different systems compatible with the measurements,
Figure~\ref{fig:plots systems from origin} shows the trajectories
emanating from the origin for a number of such systems.

\begin{figure}[ht]
  \begin{center}
    \includegraphics[width=.5\linewidth,trim=3.5cm 1cm 3.5cm 1cm, clip]{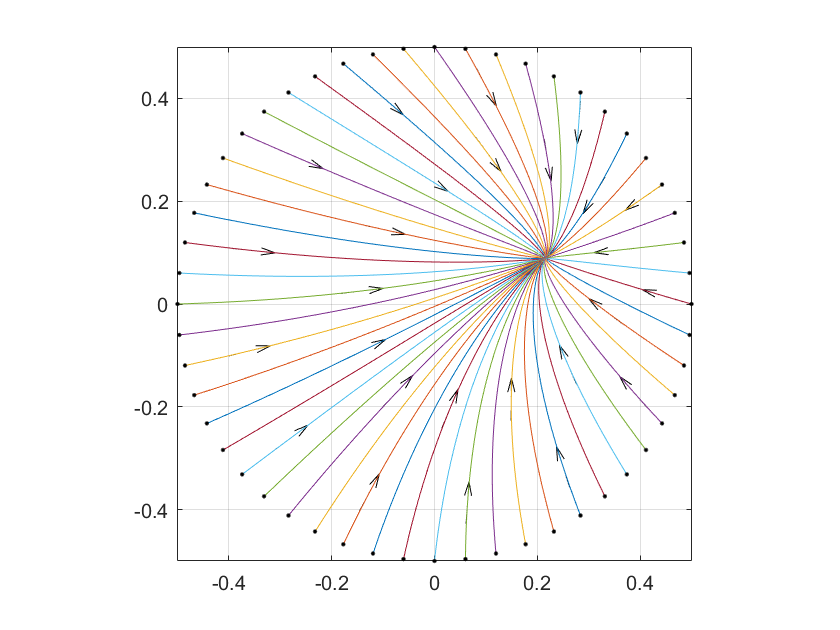}%
    \includegraphics[width=.5\linewidth,trim=3.5cm 1cm 3.5cm 1cm, clip]{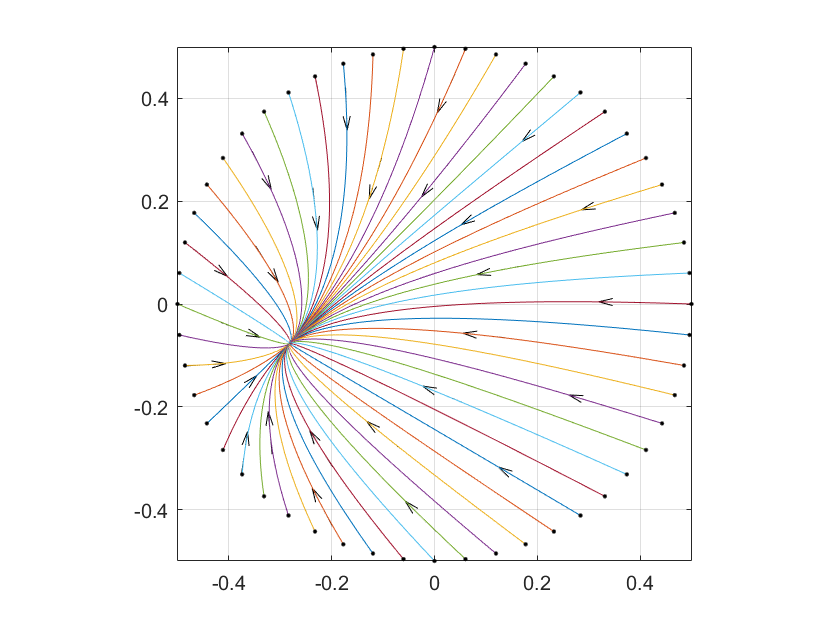}\\
    \includegraphics[width=.5\linewidth,trim=3.5cm 1cm 3.5cm 1cm, clip]{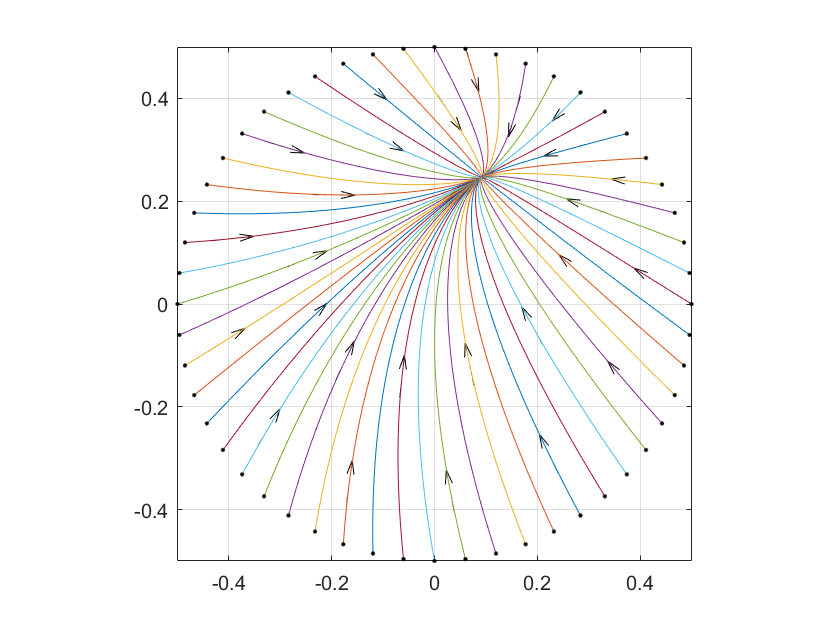}%
    \includegraphics[width=.5\linewidth,trim=3.5cm 1cm 3.5cm 1cm, clip]{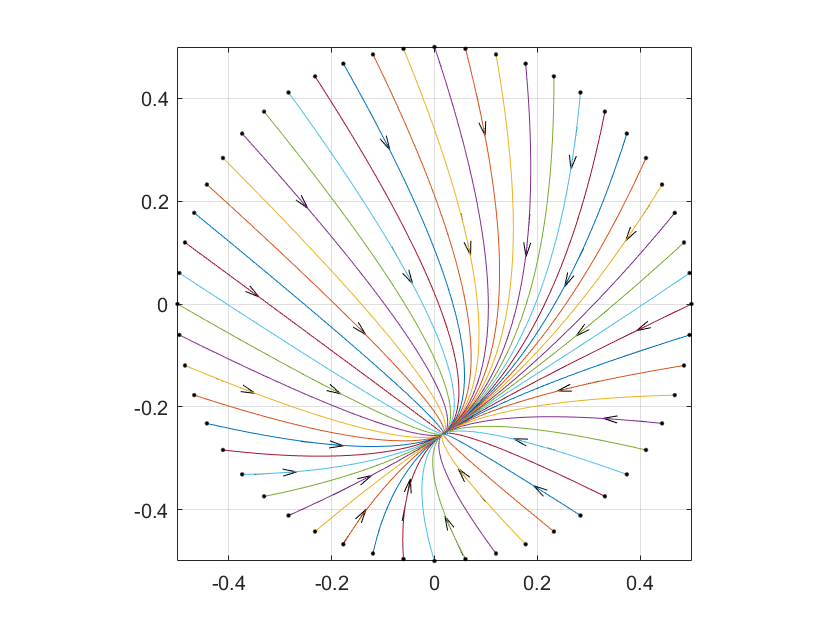}%
  \end{center}
  \caption{Trajectories of $\dot{z}=\phi^\theta(z)$ starting from a
    number of different initial conditions on the circle with radius
    $\tfrac{1}{2}$ for four different realizations of
    $\theta\in\Theta$.
  }
  \label{fig:plots systems}
\end{figure}

In order to guarantee strict contraction in terms of the measurements,
we employ the condition of Corollary~\ref{cor:contraction}. With these
measurements $\Phi$ has full row rank and thus $\Theta=\calZ(N)$ is
compact. For the last assumption, from the fact that
$\norm{J(b)}_2=\sqrt{2}$, we obtain that for any
$z,z^\star\in\mathbb{R}^2$,
\vspace{-0.3em}\[ \vspace{-0.3em}
  \norm{b(z)-b(z^\star)}_{2,(-N_{22})^{-1/2}} \leq
  \frac{\sqrt{2}}{\sqrt{-\lambda_{\textrm{max}}(N_{22})}}
  \norm{z-z^\star}_2.
\]
This means that we can use Corollary~\ref{cor:contraction}, which
shows that the unknown system is strictly contracting if
\vspace{-0.3em}\[
\vspace{-0.3em}  \osL(\phi^{\lse}) < -
  \frac{\sqrt{2\lambda_{\textup{max}}(N|N_{22})}}{\sqrt{-\lambda_{\textrm{max}}(N_{22})}}
  \approx -1.78.
\]
Using a line of reasoning as in \eqref{eq:osL phi hat} allows us to
conclude that this indeed holds. As such, the system
$\dot{z}=\phi^\theta$ is strictly contracting for all
$\theta\in\Theta$. We plot a number of trajectories of the systems
corresponding to 4 different realizations of the parameter $\theta$ in
Figure~\ref{fig:plots systems}.  It can be seen that each is strictly
contracting.

\subsection{Suboptimal regulation of unmanned aerial vehicle}
This example considers the suboptimal regulation of a fixed-wing
unmanned aerial vehicle (UAV). We consider a model of the longitudinal
motion for which the relevant equations of motion are derived in
e.g.,~\cite[Ch. 2]{JB:91}. For the values of the parameters and the
linearization we follow~\cite{EB-PT-AS:22}, which derives a benchmark
model on the basis of 10 different models of real-world UAV's of
different specifications. The derived system is a state-space model
\[
  \dot{x}(t) = A_{\textrm{cont}}x(t)+B_{\textrm{cont}}u(t),
\] 
given in~\cite[Eq. (8)-(10)]{EB-PT-AS:22}, with states
$x(t)\in\mathbb{R}^4$ and inputs $u(t)\in\mathbb{R}$. Here, $x_1(t)$
denotes the \textit{deviation from the nominal forward velocity},
$x_2(t)$ the \textit{deviation from the nominal angle of attack},
$x_3(t)$ the \textit{pitch angle}, and $x_4(t)$ the \textit{pitch
  rate}. The input $u(t)$ represents the \textit{elevator
  deflection}. We use the parameters derived in~\cite[Eq. (13)-(14)]{EB-PT-AS:22}, and thus
\vspace{-0.3em}\[
  A_{\textrm{cont}}\!\! =\!\!
  \scriptsize{
    \begin{bmatrix}
      -0.240 & 0.345
      &\!\! -0.411 & 0
      \\
      -1.905 & -10.695\!\! & 0 & \! 0.941
      \\
      0&0&0 & 1
      \\
      0.457 & -250.513 \!\! & 0
                &\! -8.844
    \end{bmatrix}
  } \!, B_{\textrm{cont}} \!\! =\!\!
  \begin{bmatrix}
    0 \\ -0.301 \\ 0 \\
    -98.658
  \end{bmatrix}
  \!.\vspace{-0.3em}
\]
We consider a discretization of this model with stepsize
$\tfrac{1}{20}${\color{black}s}, thus arriving at a system of the form
\eqref{eq:dist lin system}, where
$\hat{A}= I_4 +\tfrac{1}{20} A_{\textrm{cont}}$ and
$\hat{B}=\tfrac{1}{20} B_{\textrm{cont}}$. These matrices are unknown
to us, and instead we have access to a single second-long
\footnote{{\color{black}Note that, in the simulation we applied a fixed
    input. In reality, one would not want the UAV to be uncontrolled
    for a full second.  It should be noted that the results are also
    valid for data collected in shorter bursts.
}}
trajectory of measurements of the state $x(t)$ with
$x(0) = \tiny{\begin{pmatrix} 1&1&1&1\end{pmatrix}^\top}$, and fixed
input $u(t)=-4$. Thus, $T=20$. Denoting the matrix collecting the
noise samples by $W_s$, we assume that the noise on the measurements
is such that $W_sW_s^\top\leq 10^{-4} I_4$. These measurements give
rise to a set of systems consistent with the data given by
$\Sigma =\calZ(M)$.

Following Section~\ref{sec:subopt regulation}, we aim at finding
suboptimal values of the norm of an unknown cost function. In this
example, we consider a simple cost function
$\hat{\phi}:\mathbb{R}^4\rightarrow\mathbb{R}^4$ of the form
$\hat{\phi}(x) = x-\hat{x}$, where
$\hat{x} = \tiny{\begin{pmatrix} 1&0&0&0\end{pmatrix}^\top}$. This
situation might arise, for instance, when the UAV is tasked to mirror
the orientation and speed of another UAV on the basis of noisy
measurements. If the cost function were known, the problem of
minimizing $\norm{\hat{\phi}}_2$ would be to simply regulate $x$ as
close to $\hat{x}$ as possible. However, we assume that we only know
that $\hat{\phi}$ is affine in $x$, and thus of the form
\[
\phi^\theta(x) := \theta^\top
\begin{pmatrix}
  1
  \\
  x
\end{pmatrix}, \textrm{ where }
  \hat{\theta} =
  \begin{bmatrix}
    -\hat{x}^\top
    \\
    I_4
  \end{bmatrix}.
\]
We assume that the measurements available to us are collected
concurrently with the collection of the measurements of the unknown
system. In line with Section~\ref{sec:problem}, this means that we
measure at $z_i:=x(i-1)$ for $i=1,\ldots,T=20$. We assume a noise model
of the form $WW^\top\leq I_4$, and obtain a set of parameters
consistent with the measurements given by $\Theta=\calZ(N)$.

We are interested in finding upper bounds for the quantity
\vspace{-0.3em}\[
\vspace{-0.3em}  \min_u \norm{\hat{\phi}((I-\hat{A})^{-1}\hat{B}u)}_2.
\]
Before  applying Theorem~\ref{thm:subopt reg}, we need  to
check its assumptions. First, we employ Lemma~\ref{lem:stab} in order
to check whether each $A$ consistent with the measurements is
stable. Solving the LMI~\eqref{eq:LMIstab} yields
\vspace{-0.3em}\[
\vspace{-0.3em}  \bar{P} \approx \scriptsize{
    \begin{bmatrix} 
      \vphantom{-}0.0072 & 0.0015 & \vphantom{-}0.0004 & -0.0102 \\  
      \vphantom{-}0.0015 & 0.1665 & \vphantom{-}0.1218 & \vphantom{-}0.1496 \\
      \vphantom{-}0.0004 & 0.1218 & \vphantom{-}0.1648 & -1.3489 \\  
      -0.0102 & 0.1496 & -1.3489 			 & 45.00437  
    \end{bmatrix}}.
\] 
Next, we employ Lemma~\ref{lem:lipschitz} to obtain a Lipschitz
constant. Given that the basis functions are linear, we can easily
minimize $L$ over \eqref{eq:LMI Lip J} to obtain $L=2.0171$. Thus, $L$ is a
Lipschitz constant for $\phi^\theta$ for all $\theta\in\calZ(N)$.

Therefore, the required pieces are in place, and we can employ
Theorem~\ref{thm:subopt reg} to show that
\begin{equation}\label{eq:sims value regulation}
  \min_u
  \norm{\hat{\phi}((I-\hat{A})^{-1}\hat{B}u)}_2 \leq 2.3347.
\end{equation}
{\color{black} This bound corresponds to \eqref{eq:subopt reg} with 
minimizers $u^\star=0.2321$ and $x^\star =
(\begin{smallmatrix}
   0.809
   &0.054 & -0.355
   &-0.006
 \end{smallmatrix})^{\!\top}$ on the left- and right-hand sides, respectively.}

\begin{figure}[ht]
  \begin{center}
    \includegraphics[width=.5\linewidth,trim=1.5cm 1cm 1.5cm 1cm,
    clip]{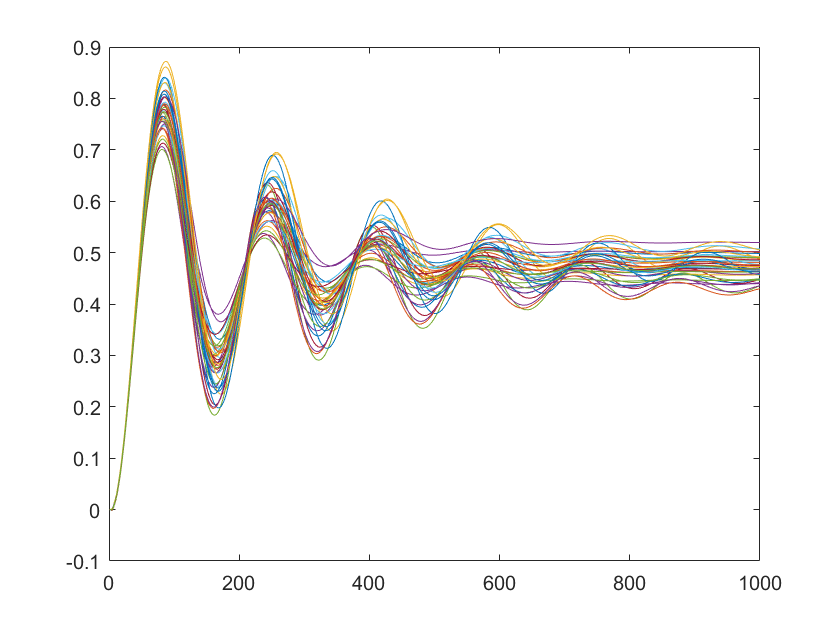}%
    \includegraphics[width=.5\linewidth,trim=1.5cm 1cm 1.5cm 1cm, clip]{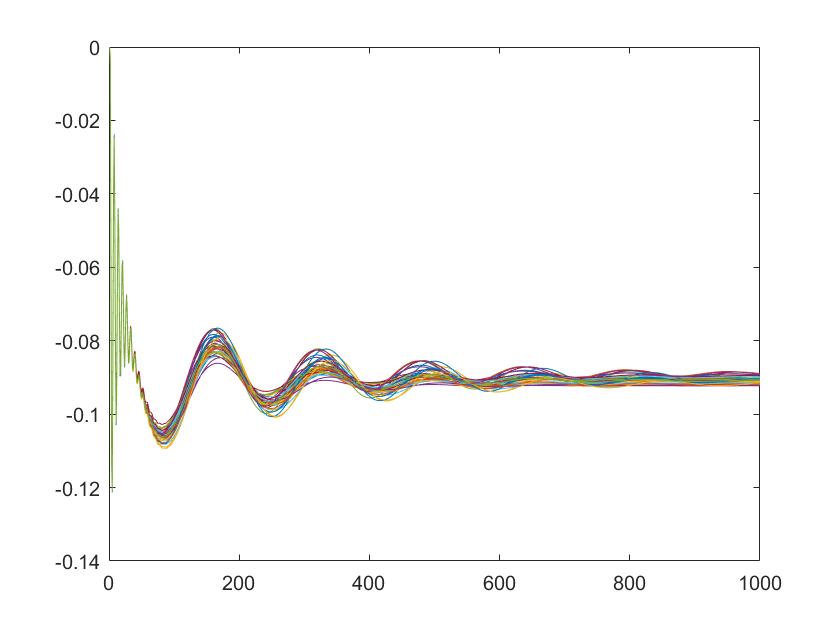}\\
    \includegraphics[width=.5\linewidth,trim=1.5cm 1cm 1.5cm 1cm,
    clip]{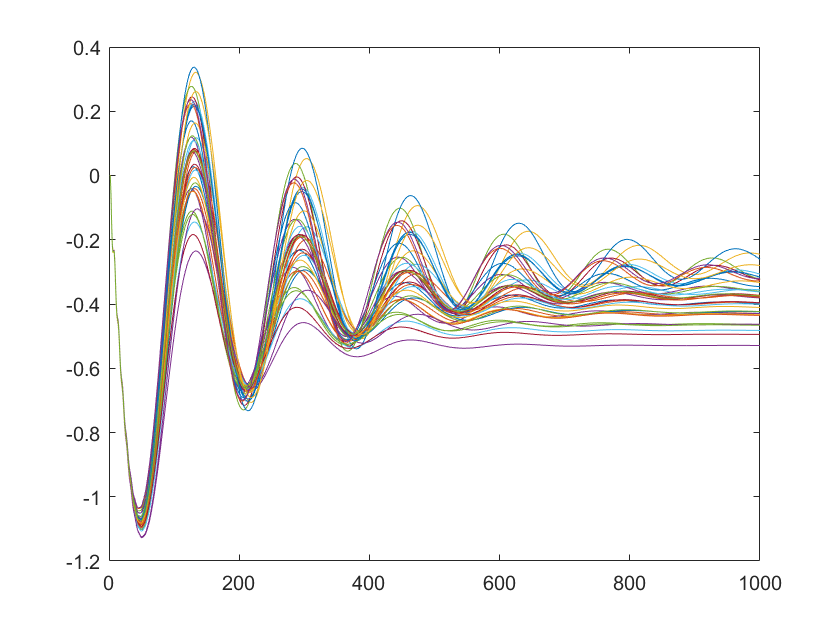}%
    \includegraphics[width=.5\linewidth,trim=1.5cm 1cm 1.5cm 1cm,
    clip]{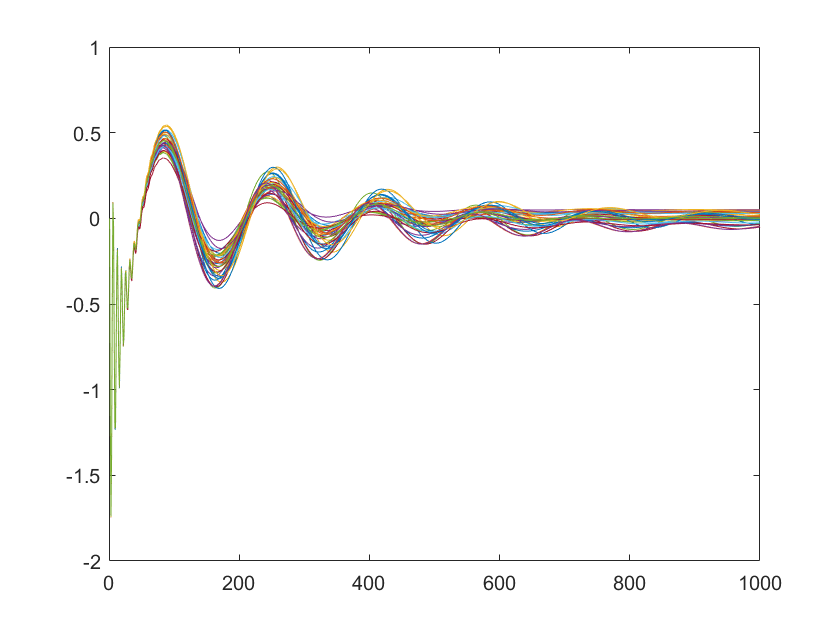}%
  \end{center}
  \caption{State trajectory corresponding to the found input
    $u^\star=0.2321$ for 40 different realizations of the matrices
    $(A,B)$ consistent with the measurements.}
  \label{fig:plots states regulation}
\end{figure}

Figure~\ref{fig:plots states regulation} shows the results of applying
this input to different realizations of the system matrices compatible with
the measurements. We make the following observations. First, one can
see that each of the systems converges to a fixed point, and each of
these fixed points are within a distance $\varepsilon^-(x)=0.4491$ of
$x$. Second, for certain systems the transients, especially in states
$x_1$ and $x_2$, converge quite slowly. This holds because
the true system is marginally stable and the same holds for elements
of set $\calZ(M)$ (to be precise, the upper bound of
Theorem~\ref{thm:transients} can only be shown to hold for
$0.999<\lambda <1$).

Corresponding to these trajectories, Figure~\ref{fig:plots norms
  regulation} shows the resulting values of
$\norm{\phi^\theta(x_k)}_2$, for both the true value and for different
realizations of $\theta\in\calZ(N)$.  The plot reveals that,
with the true parameter $\hat{\theta}$, the values at the fixed points
are well below the bound~\eqref{eq:sims value regulation}. In fact,
resolving the problem with known $\hat{\phi}$ yields
$u^\star = 0.1820$ and an upper bound of $0.7963$. Indeed, as the left
plot shows, the unknown nature of $\hat{\phi}$ has a large effect and
therefore has to be taken into account when providing formal
guarantees.

\begin{figure}[ht]
  \begin{center}
    \includegraphics[width=.5\linewidth,trim=1.5cm 1cm 1.5cm 1cm,
    clip]{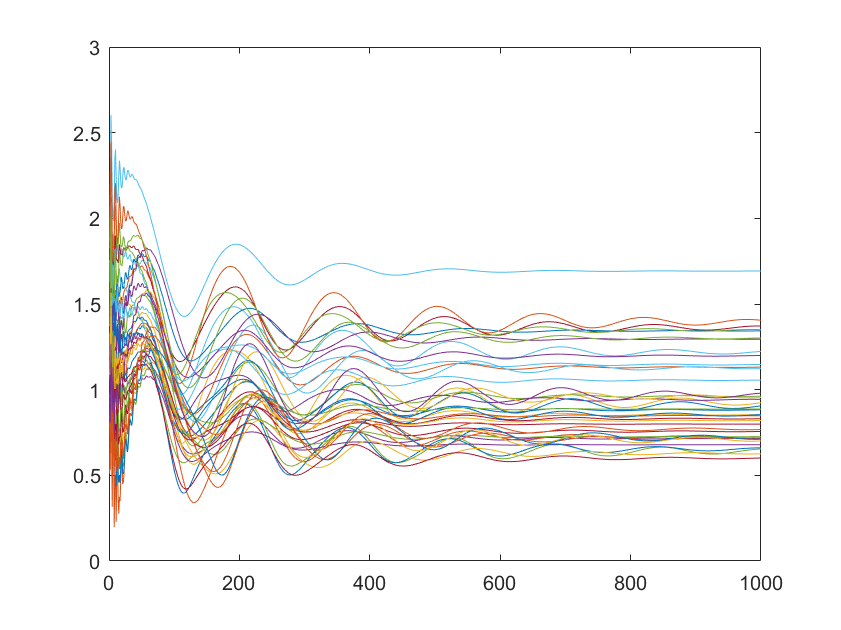}%
    \includegraphics[width=.5\linewidth,trim=1.5cm 1cm 1.5cm 1cm,
    clip]{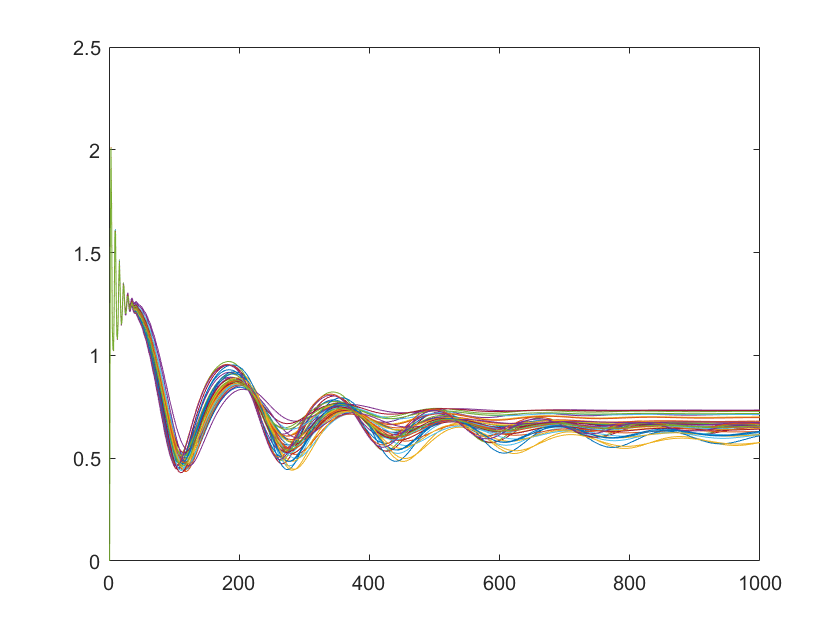}
  \end{center}
  \caption{Values of $\norm{\phi^\theta(x_k)}_2$ for the state
    trajectories of Figure~\ref{fig:plots states regulation}. On the
    left, we pick for each system a different realization of
    $\theta\in\calZ(N)$. On the right, we use the true value of
    $\hat{\theta}$.}
  \label{fig:plots norms regulation}
\end{figure}

\section{Conclusions}\label{sec:conclusions}

We have developed a set-valued regression approach to data-based
optimization of an unknown vector-valued function.  Taking a
worst-case perspective, we have provided a range of guarantees in
terms of measurements on the minimization of least conservative upper
bounds of both the norm and any linear combination of the components
of the unknown function.  Our analysis has yielded closed-form
expressions for the proposed upper bounds, conditions to ensure their
convexity (as an enabler for the use of gradient-based methods for
their optimization), and Lipschitzness characterizations to facilitate
simplifications based on interpolation.

We have illustrated the applicability of the proposed cautious
optimization approach in systems and control scenarios with unknown
dynamics: first, by providing conditions on the data which guarantee
that the dynamics is strongly contracting and, second, by providing
conditions under which we can regulate an unknown system to a point
with a guaranteed maximal value of an unknown cost function.  Finally,
we have considered online scenarios, where one repeatedly combines the
optimization of the unknown function with the collection of additional
data.  Under mild assumptions, we show that repeated measurements lead
to a strictly shrinking set of compatible parameters and we build on
this result to design an online procedure that provides sequential
upper bounds which converge to a true optimizer.

{\color{black} Future work will investigate conditions to ensure other
  regularity properties of the unknown function and its least
  conservative upper bounds. Moreover, we envision an investigation of
  the impact of the choice of basis functions on the guarantees and
  computational efficiency of the proposed cautious optimization
  methods (particularly, the use of polynomial bases coupled with
  sum-of-squares optimization techniques). The links between choices
  of basis and the performance on out-of-basis functions under certain
  regularity conditions is an interesting direction. Lastly, many
  applications in the field of control remain open, such as the
  certification of Lyapunov stability and controller design for
  nonlinear systems.}

{ \small
  \bibliographystyle{plainnat}
  \bibliography{../bib/alias,../bib/JC,../bib/Main,../bib/Main-add,../bib/New,../bib/FB}
}

\appendix

{\color{black}
  \section{Determining regularity properties using
    data}\label{sec:appendix}
\subsection{Nonnegativity of the parameters} 
The following result provides a test in terms of the data for $\theta c$ to be elementwise nonnegative for each $\theta$ compatible with the data. 
\begin{lemma}\longthmtitle{Nonnegativity of
		parameters}\label{lem:theta positive}
	Given data $(Y,\Phi)$, let $\Theta=\calZ(N)$, with $N$ is as in
	\eqref{eq:defN}. Let $c\neq 0$.  Then,
	$\theta c\in\mathbb{R}^k_{\geq 0}$ for all $\theta\in\Theta$ if and
	only if $\Phi$ has full row rank and one of the following conditions
	hold
	\begin{enumerate}
		\item\label{lemitem:single case} $c^\top(N|N_{22})c=0$ and
		$-N_{22}^{-1}N_{21}c\in\mathbb{R}^k_{\geq 0}$, or
		\item\label{lemitem:scalar} $c^\top(N|N_{22})c>0$,
		$-N_{22}^{-1}N_{21}c \in \mathbb{R}^k_{>0}$ and for all
		$i=1,\ldots,k$,
		\vspace{-0.3em}    \[
		\vspace{-0.3em}      c^\top(N|N_{22})c+\frac{(e_i^\top N_{22}^{-1}N_{21}
			c)^2}{e_i^\top N_{22}^{-1} e_i}\leq  0.
		\]
	\end{enumerate}
\end{lemma}
\vspace{-1em}\begin{pf} Recall the definition of $N_c$ from
	\eqref{eq:def Nc}. Since $c\neq 0$ we have that
	$\calZ(N)c=\calZ(N_c)$. Thus, $\theta c\in\mathbb{R}^k_{\geq 0}$ for
	all $\theta\in\Theta$ if and only if
	$\calZ(N_c)\subseteq \mathbb{R}^k_{\geq 0}$.
	
($\Rightarrow$): Suppose that
		$\calZ(N_c)\subseteq \mathbb{R}^k_{\geq 0}$.  This implies that
		the set $\calZ(N_c)$ contains no {\color{black} nontrivial}
		subspace, which implies that $N_{22}<0$. In turn, this holds only
		if $\Phi$ has full row rank.
		
		Recall that $N|N_{22}\geq 0$ by assumption. Therefore we have
		either $c^\top(N|N_{22})c=0$ or $c^\top(N|N_{22})c>0$.  If
		$c^\top(N|N_{22})c=0$, then $N_c\leq 0$. As such,
		$\calZ(N_c)=\{ -N_{22}^{-1}N_{21}\}$, that is, \textit{only} the
		least-squares estimate is compatible with the data. Hence,
		\emph{\ref{lemitem:single case}} follows.
		
		Next, consider $c^\top(N|N_{22})c>0$, that is, $N_c$ has one
		positive eigenvalue. Note that the Slater condition holds, and we
		can apply the S-Lemma~\cite[Thm 2.2]{IP-TT:07} to see that
		$\calZ(N_c)\subseteq \mathbb{R}^k_{\geq 0}$ is equivalent to the
		existence of $\alpha_1,\ldots,\alpha_k\geq 0$ with
		\vspace{-0.3em}
		\[
		\vspace{-0.3em}
		\begin{bmatrix}
			0 & e_i^\top
			\\
			e_i & 0
		\end{bmatrix}
		-\alpha_iN_c\geq 0,
		\]
		for each standard basis vector $e_i$.  This requires
		$\alpha_i\neq 0$ for each $i \in \until{k}$. Thus, this can
		equivalently written as
		\[
		N_c -\frac{1}{\alpha_i}
		\begin{bmatrix}
			0 & e_i^\top
			\\
			e_i &
			0
		\end{bmatrix} \leq 0.
		\]
		Since $N_{22}<0$, we can define $\beta_i=\tfrac{1}{\alpha_i}$ and
		conclude that the latter holds if and only if
		\begin{align*}
			& c^\top N_{11} c - (c^\top N_{12} -\beta_i e_i^\top) N_{22}^{-1}
			(N_{21} c - \beta_i e_i)
			\\
			& \quad = c^\top(N|N_{22})c + \beta_i (c^\top N_{12} N_{22}^{-1}
			e_i + e_i^\top N_{22}^{-1}
			N_{21} c)
			\\
			& \qquad - \beta_i^2 e_i^\top N_{22}^{-1} e_i
			\leq 0 .    
		\end{align*}
		In turn, there exist such $\beta_i$ if and only if the minimal
		value over all $\beta_i>0$ satisfies this. Since this is a scalar
		expression, we can explicitly minimize it to prove that
		\emph{\ref{lemitem:scalar}} holds.
		
		($\Leftarrow$): Similar to before, if $\Phi$ has full row rank and
		\emph{\ref{lemitem:single case}} holds, then
		$\calZ(N_c)=\{ -N_{22}^{-1}N_{21}\}$. Hence
		$\calZ(N_c)\subseteq \mathbb{R}^k_{\geq 0}$ follows
		immediately. If, instead, $\Phi$ has full row rank and
		\emph{\ref{lemitem:scalar}} holds, reversing the steps of the
		previous part of the proof yields that
		$\calZ(N_c)\subseteq \mathbb{R}^k_{\geq 0}$. \QED 
\end{pf}

  \subsection{Bounding the Jacobian}
  Let $f:\mathbb{R}^n\rightarrow\mathbb{R}^m$ be differentiable and
  denote its Jacobian,
  \[
    J(f) :=
    \begin{bmatrix} \tfrac{df}{dz_1}
      & \cdots &
                 \tfrac{df}{dz_n}
    \end{bmatrix}
    =
    \scriptsize
    {
      \begin{bmatrix}
        \tfrac{df_1}{dz_1}
        & \cdots & \tfrac{df_1}{dz_n}
        \\ \vdots & & \vdots
        \\
        \tfrac{df_m}{dz_1}
        & \cdots &
                   \tfrac{df_m}{dz_n}
      \end{bmatrix}}.
  \]
  Clearly, $ J(\phi^\theta) = \theta^\top J(b)$.  The following
  result, consequence of Theorem~\ref{thm:explicit delta}, provides an
  expression for the Jacobian of the worst-case linear bound function.
  
  \begin{corollary}\longthmtitle{Jacobians of the worst-case linear 
      bound}\label{cor:gradient}
    Given data $(Y,\Phi)$, let $\Theta=\calZ(N)$ and assume $N_{22}<0$.
    Let the basis functions $\phi_i(z)$ be differentiable and such that
    $b(z) \neq0$ for all $z$.  Then
    \[\begin{split}
        &J(g_c(z;\Theta)) =
        \\
        &\bigg(\!c^\top N_{12}
          + \frac{\sqrt{c^\top
          (N|N_{22})c}}{\sqrt{b(z)^{\!\top} (-N_{22}^{-1})
          b(z)}}b(z)^{\!\top}\bigg)(-N_{22}^{-1})J(b(z)).    
      \end{split}\]
  \end{corollary}	
  In a similar fashion, we can also obtain gradients of the bound of
  $g(z;\Theta)$ given in \eqref{eq:expl g}. Note that, if $N_{22}<0$
  this requires that both $b(z)\neq 0$ and
  $\phi^{\lse}(z;\Theta)\neq 0$ for all $z$.

\subsection{Lipschitz continuity}

We investigate Lipschitz continuity for an unknown
function~$\hat{\phi}$.
In line with the rest of the paper, we guarantee this by providing
conditions for the function $\phi^\theta$ for \emph{all}
parameters~$\theta\in\Theta$.  Without added difficulty, the following
result considers the slightly more general case of weighted norms.

\begin{lemma}\longthmtitle{Lipschitz constants of functions consistent
    with the measurements}\label{lem:lipschitz}
  Given data $(Y,\Phi)$, let $\Theta=\calZ(N)$, with $N$ as
  in~\eqref{eq:defN}, and assume $N$ has at least one positive
  eigenvalue. Let $P\in\mathbb{S}^{m}$ and
  $Q\in\mathbb{S}^{n}$ with $P>0$ and $Q>0$. For $L\geq 0$ and
  $z,z^\star\in\calS\subseteq\mathbb{R}^n$, we have
  \[
    \norm{\phi^\theta(z)-\phi^\theta(z^\star)}_{2,P^{1/2}} \leq L
    \norm{z-z^\star}_{2,Q^{1/2}} \quad \textrm{ for all } \quad
    \theta\in\Theta
  \]
  if and only if there exists $\alpha\geq 0$ such that
  \vspace{-0.3em}\begin{equation}\label{eq:Lipschitz}
    \scriptsize{\begin{pmat}[{.|}]
        \! L^2 P^{-1} \! &0 &0 \cr
        0& 0 &\frac{b(z)-b(z^\star)}{\norm{z-z^\star}_{2,Q^{1/2}}} \cr\-
        0& \frac{b(z)^\top-b(z^\star)^\top}{\norm{z-z^\star}_{2,Q^{1/2}}} & 1 \cr
      \end{pmat}} \!\! - \alpha
    \normalsize{\begin{pmat}[{|}]
        N & 0 \cr\- 0 & 0 \cr
      \end{pmat}}\geq 0. \vspace{-0.3em} 
  \end{equation}
  If, in addition the basis functions $\phi_i$ are differentiable, then 
  \vspace{-0.2em}\[\vspace{-0.2em}
    \norm{J(\phi^\theta)}_2 \leq L \quad \textrm{ for all } \quad
    \theta\in\Theta
  \]	
  if and only if there exists $\alpha\geq 0$ such that
  \vspace{-0.2em}\begin{equation}\label{eq:LMI Lip J}
    \vspace{-0.2em}\scriptsize{\begin{pmat}[{.|}]
        L^2 I_m &0 &0 \cr
        0& 0 &J(b) \cr\-
        0& J(b)^\top & I_n \cr
      \end{pmat} - \normalsize{\alpha
        \begin{pmat}[{|}]
          N & 0 \cr\- 0 & 0 \cr
        \end{pmat}}\geq 0.}
  \end{equation}
\end{lemma}
\begin{pf}
  Note that
  $\phi^\theta(z)-\phi^\theta(z^\star) =
  \theta^\top(b(z)-b(z^\star))$. Similarly,
  $J(\phi^\theta)= \theta^\top J(b)$. The result can then be proven by
  following the steps of the proof of Lemma~\ref{lem:lmi for
    subopt}. \QED
\end{pf}

\begin{remark}\longthmtitle{Single check for establishing Lipschitz
    constant}\label{rem:check for Lip}
  {\rm To check for Lipschitzness, we can employ the method of
    Remark~\ref{rem:overestimating} to avoid checking
    \eqref{eq:Lipschitz} for all $z,z^\star\in\mathbb{R}^n$. Suppose
    the basis satisfies a Lipschitz condition of the form
    \[ \norm{b(z)-b(z^\star)}_2 \leq L_b\norm{z-z^\star}^2_{2,Q^{1/2}},
    \]
    for all $z,z^\star\in\mathbb{R}^n$.  Then,
    \[
      \vspace{-0.2em}
      \frac{(b(z)-b(z^\star))(b(z)-b(z^\star))^\top}{\norm{z-z^\star}^2_{2,Q^{1/2}}}
      \leq L_b^2I,\vspace{-0.2em}
    \]
    and thus, as in Remark~\ref{rem:overestimating}  we have
    \begin{equation}\label{eq:LMI for Lip}
      \scriptsize{
        \begin{bmatrix}
          L^2P^{-1} & 0
          \\
          0&
             -L_b^2I_k \end{bmatrix}} - \normalsize{\alpha N\geq 0}
       \end{equation}
       implies \eqref{eq:Lipschitz}. This gives a single condition to
       check to establish a Lipschitz constant valid for all
       $z,z^\star\in\mathbb{R}^n$.  \oprocend }
   \end{remark}
 }

\end{document}